\documentclass[10pt,a4paper]{amsart}
\usepackage{amssymb,amsmath}
\usepackage[american,french]{babel}
\usepackage{amsopn}
\usepackage{graphicx}
\usepackage{fancyhdr}
\usepackage[T1]{fontenc}
\title{About quadratic differential operators}
\newcommand{\rr}{\mathbb{R}}
\newcommand{\eps}{\varepsilon}
\newcommand{\nn}{\mathbb{N}}
\newcommand{\cc}{\mathbb{C}}

\def\un{{\mathrm{1~\hspace{-1.4ex}l}}}
\def\init{\setcounter{equa}{0}}
\def\inc{\stepcounter{equa}}
\def\num{\tag{\thesection.\theequa}}

\def\wrtext#1{\relax\ifmmode{\leavevmode\hbox{#1}}\else{#1}\fi}
\def\abs#1{\left|#1\right|}
\def\begeq{\begin{equation}}
\def\endeq{\end{equation}}

\def\part#1{\frac{\partial}{\partial #1}}

\def\norm#1{||\,#1\,||}

\newcommand{\real}{\mathbb{R}}
\newcommand{\comp}{\mathbb{C}}

\newcommand{\nat}{\mathbb{N}}

\renewcommand{\Re}{\textrm{Re }}
\renewcommand{\Im}{\textrm{Im }}
\renewcommand{\exp}{\mbox{\rm exp\,}}

\def\neigh{neighborhood}

\def\Re{{\rm Re\,}}
\def\Im{{\rm Im\,}}

\newtheorem{dref}{Definition}
\newtheorem{theo}[dref]{Theorem}
\newtheorem{prop}[dref]{Proposition}

\begin{document}
\newcounter{equa}
\selectlanguage{american}

\newtheorem{lemma}{Lemma}[subsection]
\newtheorem{definition}{Definition}[subsection]
\newtheorem{proposition}{Proposition}[subsection]
\newtheorem{theorem}{Theorem}[subsection]

\title[SEMICLASSICAL ESTIMATES FOR OPERATORS WITH DOUBLE CHARACTERISTICS]{SEMICLASSICAL HYPOELLIPTIC ESTIMATES FOR NON-SELFADJOINT OPERATORS WITH DOUBLE CHARACTERISTICS}
\author{Michael Hitrik, Karel Pravda-Starov}
\address{\noindent UCLA Department of Mathematics, Los Angeles, CA 90095-1555, USA}
\email{hitrik@math.ucla.edu}
\urladdr{http://www.math.ucla.edu/~hitrik/}

\address{Department of Mathematics,
Imperial College London,
Huxley Building, 180 Queen's Gate,
London SW7 2AZ, UK}
\email{k.pravda-starov@imperial.ac.uk}
\urladdr{http://www.ma.ic.ac.uk/~kpravdas/}

\begin{abstract}
For a class of non-selfadjoint semiclassical pseudodifferential operators with double characteristics, we study bounds for resolvents and estimates for low lying eigenvalues. Specifically, assuming that the quadratic approximations of the principal symbol of the operator along the double characteristics enjoy a partial ellipticity property along a suitable subspace of the phase space, namely their singular spaces, we establish
semiclassical hypoelliptic a priori estimates with a loss of the full power of the semiclassical parameter giving a localization for the low lying spectral values of the operator.
\end{abstract}


\maketitle

\section{Introduction}
\init

\subsection{Miscellaneous facts about quadratic differential operators and doubly characteristic pseudodifferential operators}

Since the classical work by J. Sj\"ostrand~\cite{sjostrand}, the study of spectral properties of quadratic diffe\-rential
operators has played a basic r\^ole in the analysis of partial differential operators with double characteristics.
Roughly speaking, if we have, say, a classical pseudodifferential operator $p^w(x,D_x)$ on $\rr^n$ with the Weyl
symbol
$p(x,\xi)=p_m(x,\xi)+p_{m-1}(x,\xi)+\ldots$ of order $m$, and if $X_0=(x_0,\xi_0)\in \rr^{2n}$ is a point
where
$$p_m(X_0)=dp_m(X_0)=0,$$
then it is natural to consider the quadratic form $q$ which begins the Taylor expansion of
$p_m$ at $X_0$ in order to investigate the properties of the pseudodifferential operator $p^w(x,D_x)$. For example, the study of a priori estimates such as hypoelliptic estimates of the form
$$
\norm{u}_{m-1}\leq C_K \left(\norm{p^w(x,D_x) u}_0+\norm{u}_{m-2}\right),\quad u\in C^{\infty}_0(K), \quad
K\subset \subset \rr^n,
$$
then often depends on the spectral analysis of the quadratic operator $q(x,\xi)^w$.
See~\cite{hypoelliptic}, as well as Chapter 22 of~\cite{hormander} together with further references given there.
In the classical work~\cite{sjostrand}, the spectrum of a general quadratic differential operator, that is an operator defined in the Weyl quantization
\begin{equation}\label{3}\inc
q^w(x,D_x) u(x) =\frac{1}{(2\pi)^n}\int_{\rr^{2n}}{e^{i(x-y).\xi}q\Big(\frac{x+y}{2},\xi\Big)u(y)dyd\xi}, \num
\end{equation}
by a symbol $q(x,\xi)$, where $(x,\xi) \in \rr^{n} \times \rr^n$ and $n \in \nn^*$, which
is a complex-valued quadratic form, has been determined under the
basic assumption of global ellipticity of the quadratic symbol
$$(x,\xi) \in \rr^{2n}, \ q(x,\xi)=0 \Rightarrow (x,\xi)=0.$$

\medskip

We recently investigated properties of non-elliptic quadratic operators in the works \cite{HiPr} and~\cite{karel}. Considering quadratic operators whose Weyl symbols have real parts with a sign, say here, Weyl symbols with non-negative real parts
\begin{equation}\label{smm1}\inc
\textrm{Re }q \geq 0, \num
\end{equation}
we pointed out the existence of a particular linear subvector space $S$ in the phase space $\rr_x^n \times \rr_{\xi}^n$ intrinsically associated to their Weyl symbols $q(x,\xi)$ and
called singular space, which seems to play a basic r\^ole in the understanding of a number of fairly general properties, such as spectral or subelliptic properties, of these non-elliptic quadratic operators.
In particular, we established that when a quadratic symbol fulfilling (\ref{smm1}) satisfies an assumption of partial ellipticity along its singular space $S$, that is,
$$(x,\xi) \in S, \ q(x,\xi)=0 \Rightarrow (x,\xi)=0,$$
then the spectrum of the quadratic operator $q^w(x,D_x)$ is only composed of a countable number of eigenvalues of finite multiplicity, with a structure
similar to the one known in the case of global ellipticity described by J.~Sj\"ostrand in~\cite{sjostrand}.

\medskip

The purpose of the present work is to address the question of how these recent improvements in the understanding of non-elliptic quadratic operators allow to enhance the comprehension of the properties of certain classes of non-selfadjoint semiclassical operators with double characteristics. In this work, which is planned to be the first one in a series on doubly characteristic pseudodifferential operators, we shall study bounds for resolvents and estimates for low lying eigenvalues for non-selfadjoint semiclassical pseudodifferential operators with principal symbols whose quadratic approximations at doubly characteristic points enjoy a partial ellipticity property along their singular spaces. Under these particular assumptions of  partial ellipticity for these quadratic approximations, we shall establish a semiclassical hypoelliptic a priori estimate with a loss of the full power of the semiclassical parameter which gives a localization for the low lying spectral values of the operator.

\medskip

Before giving the precise statement of our main result, we shall recall miscellaneous facts and notation
that we will need about quadratic differential operators. Associated to a complex-valued quadratic form
\begin{eqnarray*}
q : \rr_x^n \times \rr_{\xi}^n &\rightarrow& \cc\\
 (x,\xi) & \mapsto & q(x,\xi),
\end{eqnarray*}
with $n \in \nn^*$, is the Hamilton map $F \in M_{2n}(\cc)$ uniquely defined by the identity
\begin{equation}\label{10}\inc
q\big{(}(x,\xi);(y,\eta) \big{)}=\sigma \big{(}(x,\xi),F(y,\eta) \big{)}, \ (x,\xi) \in \rr^{2n},  (y,\eta) \in \rr^{2n}, \num
\end{equation}
where $q\big{(}\textrm{\textperiodcentered};\textrm{\textperiodcentered} \big{)}$ stands for the polarized form
associated to the quadratic form $q$ and $\sigma$ is the canonical symplectic form on $\rr^{2n}$,
\begin{equation}
\label{11}\inc
\sigma \big{(}(x,\xi),(y,\eta) \big{)}=\xi.y-x.\eta, \ (x,\xi) \in \rr^{2n},  (y,\eta) \in \rr^{2n}. \num
\end{equation}
It follows directly from the definition of the Hamilton map $F$ that
its real and imaginary parts, denoted respectively by $\textrm{Re } F$ and $\textrm{Im }F$, are the Hamilton maps associated
to the quadratic forms $\textrm{Re } q$ and $\textrm{Im }q$, respectively; and that a
Hamilton map is always skew-symmetric with respect to $\sigma$. This fact is just a consequence of the
properties of skew-symmetry of the symplectic form and symmetry of the polarized form
\begin{equation}\label{12}\inc
\forall X,Y \in \rr^{2n}, \ \sigma(X,FY)=q(X;Y)=q(Y;X)=\sigma(Y,FX)=-\sigma(FX,Y).\num
\end{equation}
We defined in~\cite{HiPr} the singular space $S$ associated to a quadratic symbol $q$ as the following intersection of kernels
\begin{equation}\label{h1}\inc
S=\Big(\bigcap_{j=0}^{2n-1}\textrm{Ker}\big[\textrm{Re }F(\textrm{Im }F)^j \big]\Big) \cap \rr^{2n}, \num
\end{equation}
where $F$ stands for its Hamilton map; and we proved (Theorem~1.2.2 in~\cite{HiPr}) that when a quadratic symbol $q$ with a non-negative real part is elliptic on its singular space $S$,
\begin{equation}\label{sm2}\inc
(x,\xi) \in S, \ q(x,\xi)=0 \Rightarrow (x,\xi)=0, \num
\end{equation}
then the spectrum of the quadratic operator  $q^w(x,D_x)$
is only composed of eigenvalues of finite multiplicity
\begin{equation}\label{sm6}\inc
\sigma\big{(}q^w(x,D_x)\big{)}=\Big\{ \sum_{\substack{\lambda \in \sigma(F), \\  -i \lambda \in \cc_+
\cup (\Sigma(q|_S) \setminus \{0\})
} }
{\big{(}r_{\lambda}+2 k_{\lambda}
\big{)}(-i\lambda) : k_{\lambda} \in \nn}
\Big\}, \num
\end{equation}
where $r_{\lambda}$ is the dimension of the space of generalized eigenvectors of $F$ in $\cc^{2n}$
belonging to the eigenvalue $\lambda \in \cc$,
$$\Sigma(q|_S)=\overline{q(S)} \textrm{ and } \cc_+=\{z \in \cc : \textrm{Re }z>0\}.$$
Let us finally end these few recollections by mentioning that one can also describe the singular spaces of such quadratic symbols (see Section~1.4 in~\cite{HiPr})
in terms of the eigenspaces associated to the real eigenvalues of their Hamilton maps. Considering such a quadratic symbol $q$, the set of real eigenvalues of its Hamilton map $F$ can then be written as
$$\sigma(F) \cap \rr =\{\lambda_1,...,\lambda_r,-\lambda_1,...,-\lambda_r\},$$
with $\lambda_j \neq 0$ and $\lambda_j \neq \pm \lambda_k$  if $j \neq k$; and one can check that its singular space is the direct sum of the symplectically orthogonal spaces
\begin{equation}\label{sm5bis}\inc
S=S_{\lambda_1} \oplus^{\sigma \perp} S_{\lambda_2} \oplus^{\sigma \perp}... \oplus^{\sigma \perp} S_{\lambda_r}, \num
\end{equation}
where the spaces $S_{\lambda_{j}}$, $1 \leq j \leq r$, are the symplectic spaces
\begin{equation}\label{sm5bis1}\inc
S_{\lambda_j}=\big(\textrm{Ker}(F -\lambda_j) \oplus \textrm{Ker}(F+\lambda_j) \big) \cap \rr^{2n}. \num
\end{equation}

\subsection{Statement of the main result}

Let us now state the main result contained in this paper.
Let $m\geq 1$  be a $C^{\infty}$ order function on $\real^{2n}$ fulfilling
\begin{equation}\label{eq1.1}\inc
\exists C_0 \geq 1, N_0>0,\ m(X)\leq C_0 \langle{X-Y\rangle}^{N_0} m(Y),\ X,Y\in \real^{2n}, \num
\end{equation}
where $\langle X \rangle=(1+|X|^2)^{\frac{1}{2}}$, and $S(m)$ be the symbol class
$$S(m)=\left\{ a\in C^{\infty}(\real^{2n},\cc): \forall \alpha \in \nn^{2n}, \exists C_{\alpha}>0, \forall X \in \rr^{2n}, \  |\partial_X^{\alpha} a(X)| \leq C_{\alpha} m(X)\right\}.$$
We shall assume in the following, as we may, that $m$ belongs to its own symbol class $m\in S(m)$.

Considering a symbol $P(x,\xi;h)$ with a semiclassical asymptotic expansion in the symbol class $S(m)$,
\begin{equation}\label{xi1}\inc
P(x,\xi;h) \sim \sum_{j=0}^{+\infty}{p_j(x,\xi)h^j}, \num
\end{equation}
with $p_j \in S(m)$, $j \in \nn$; such that its principal symbol $p_0$ has a non-negative real part
\begin{equation}\label{eq1.4}\inc
\textrm{Re }p_0(X)\geq 0,\ X=(x,\xi) \in \real^{2n}, \num
\end{equation}
we shall study the operator
\begin{equation}\label{eq1.3}\inc
P = P^w(x,hD_x;h),\ 0<h\leq 1, \num
\end{equation}
defined by the $h$-Weyl quantization of the symbol $P(x,\xi;h)$, that is, the Weyl quantization of the symbol $P(x,h\xi;h)$.  When equipped with the domain
$$H(m)=\big(m^w(x,hD_x)\big)^{-1}\big(L^2(\real^n)\big),$$
for $h>0$ sufficiently small, the operator $P$ becomes a closed and densely defined operator on $L^2(\real^n)$ (see Section~3 in~\cite{hager}).
We shall assume that the real part of the principal symbol $p_0$ is elliptic at infinity in the sense that
\begin{equation}\label{eq1.5}\inc
\exists C>0,\forall \abs{X}\geq C, \ \textrm{Re } p_0(X)\geq \frac{m(X)}{C}. \num
\end{equation}
This assumption ensures (see Section~3 in~\cite{hager}) that for sufficiently small values of the semiclassical parameter $h$, $0 < h \ll 1$, the spectrum of the operator $P$ in a fixed \neigh{} of $0\in \comp$ is discrete and consists of eigenvalues of finite algebraic multiplicity.

We shall also assume that the characteristic set of the real part of the principal symbol $p_0$,
$$(\textrm{Re }p_0)^{-1}(0)\subset \real^{2n},$$
is finite, so that we may write it as
\begin{equation}\label{eq1.6}\inc
(\textrm{Re } p_0)^{-1}(0) = \{X_1,...,X_N\}.\num
\end{equation}
The sign assumption (\ref{eq1.4}) implies in particular that we have
$$d\textrm{Re }p_0(X_j)=0,$$
for all $1 \leq j \leq N$; and we shall actually assume that these points are all doubly characteristic points for the symbol $p_0$,
\begin{equation}\label{eq1.6.5}\inc
p_0(X_j)=dp_0(X_j)=0,\ 1\leq j\leq N, \num
\end{equation}
so that we may write
\begin{equation}\label{eq1.6.6}\inc
p_0(X_j+Y)=q_j(Y)+\mathcal{O}(Y^3),\num
\end{equation}
when $Y\rightarrow 0$; where $q_j$ is the quadratic approximation which begins the Taylor expansion of the principal symbol $p_0$ at $X_j$. Notice that the sign assumption (\ref{eq1.4}) also implies that these complex-valued quadratic forms $q_j$ have non-negative real parts
\begin{equation}\label{kps1}\inc
\Re q_j\geq 0.\num
\end{equation}
By denoting $S_j$ the singular spaces associated to these quadratic forms $q_j$, the purpose of this work is to establish the following result:

\bigskip

\begin{theo}\label{theo}
Consider a symbol $P(x,\xi;h)$ with a semiclassical expansion in the class $S(m)$ such that its principal symbol $p_0$ fulfills the assumptions {\rm (\ref{eq1.4})}, {\rm (\ref{eq1.5})}, {\rm (\ref{eq1.6})} and {\rm (\ref{eq1.6.5})}. When all the quadratic forms $q_j$, $1 \leq j \leq N$, defined in {\rm (\ref{eq1.6.6})} are elliptic on their associated singular spaces
\begin{equation}\label{kps2}\inc
X \in S_j, \ q_j(X)=0 \Rightarrow X=0,\num
\end{equation}
then for any constant $C>1$ and any fixed \neigh{} $\Omega_j \subset \comp$ of the spectrum of the quadratic operator associated to the quadratic symbol~$q_j$,
$$\sigma\big(q_j^w(x,D_x)\big) \subset \Omega_j,$$
described in {\rm (\ref{sm6})}, there exist some positive constants $0 < h_0 \leq 1$ and $C_0>0$ such that for all $0<h \leq h_0$, $u \in \mathcal{S}(\real^n)$ and $\abs{z}\leq C$ satisfying
$$z-p_1(X_j) \notin \Omega_j, \ 1 \leq j \leq N,$$
we have
\begin{equation}\label{eq1.7}\inc
h\|u\|\leq C_0\|(P-hz)u\|,\num
\end{equation}
with $P=P^w(x,hD_x;h)$; where $p_1(X_j)$ stands for the value of the subprincipal symbol of the symbol $P(x,\xi;h)$ evaluated at the doubly characteristic point $X_j$ and $\|\cdot\|$ is $L^2$--norm on $\real^n$.
\end{theo}

\bigskip

Let us begin our few comments about Theorem~\ref{theo} by mentioning that its result was essentially well-known in the case when the quadratic forms $q_j$ are all globally elliptic on $\real^{2n}$,
$$X \in \rr^{2n}, \ q_j(X)=0 \Rightarrow X=0,$$
when $1 \leq j \leq N$. We refer the reader to the work \cite{sjostrand} of J. Sj\"ostrand where the case of classical pseudodifferential operators is considered. The novelty of Theorem~\ref{theo} comes therefore from the fact that the semiclassical hypoelliptic a priori estimate with a loss of the full power of the semiclassical parameter (\ref{eq1.7}) remains valid in cases where the global ellipticity of the Hessians of the principal symbol at doubly characteristic points fails. Our result actually shows that this estimate holds only under the weaker assumption of partially ellipticity (\ref{kps2}) for the Hessians of the principal symbol at doubly characteristic points. Let us also stress the fact that Theorem~\ref{theo} actually extends the result of J.~Sj\"ostrand in \cite{sjostrand} since one can check from the definitions (\ref{10}) and (\ref{h1}) that the singular space $S$ of a complex-valued quadratic form $q$ with a non-negative real part is always distinct from the whole phase space $\rr^{2n}$ as soon as its real part is a non-zero quadratic form
$$\exists (x_0,\xi_0) \in \rr^{2n}, \ \textrm{Re }q(x_0,\xi_0) \neq 0.$$
A noticeable example of non-elliptic quadratic operator fulfilling the assumption of partial ellipticity (\ref{kps2}) is given by the Kramers-Fokker-Planck operator
$$K=-\Delta_v+\frac{v^2}{4}-\frac{1}{2}+v.\partial_x-\big(\partial_xV(x)\big).\partial_v, \ (x,v) \in \rr^{2},$$
with the quadratic potential
$$V(x)=\frac{1}{2}ax^2, \ a \in \rr^*.$$
One can actually check that this operator can be expressed as
$$K=q^w(x,v,D_x,D_v)-\frac{1}{2},$$
with the Weyl symbol
\begin{equation}\label{kps3}\inc
q(x,v,\xi,\eta)=\eta^2+\frac{1}{4}v^2+i(v \xi-a x \eta),\num
\end{equation}
which is a non-elliptic complex-valued quadratic form with a non-negative real part and a zero singular space. Starting from this example, we may easily construct models for Hessians with non-negative real parts fulfilling the condition (\ref{kps2}) whose singular spaces $S$ are both
non-trivial and distinct of the whole phase space. Such a model is for instance obtained when adding to the quadratic form $q$ defined in (\ref{kps3}) an elliptic purely imaginary-valued quadratic form $i\tilde{q}$ in other symplectic variables $(x'',\xi'')$,
$$Q(x',x'',\xi',\xi'') = q(x',\xi') + i \tilde{q}(x'',\xi''),$$
since the singular space is in this case given by
$$S=\{(x',x'',\xi',\xi'') \in \rr^{2n'+2n''}: x'=\xi'=0\}.$$

\medskip

About the present work, we drew our inspiration quite exclusively from the semiclassical analysis for Kramers-Fokker-Planck equation led by F.~H\'erau, J.~Sj\"ostrand and C.~Stolk in \cite{HeSjSt}. Our proof of Theorem~\ref{theo} relies on a similar construction of a global bounded weight function $G$ with controlled derivatives and the use, on the FBI-Bargmann side, of associated weighted spaces of holomorphic functions on which the quadratic approximations at critical points of the new principal symbol of the operator
$$\tilde{p}_0 \sim p_0+i\delta H_G p_0,$$
become globally elliptic although the quadratic approximations of the original principal symbol may fail global ellipticity since they only fulfill the assumption of partial ellipticity on their singular spaces. The structure of our proof will therefore follow the one of the analysis led in \cite{HeSjSt} for the proof of the first a priori estimate in Theorem~1.2. Parts of our proof will actually be the same and we shall therefore refer directly the reader to some parts of the work by F.~H\'erau, J.~Sj\"ostrand and C.~Stolk when no change of any kind is needed. In the setting considered in \cite{HeSjSt}, the authors make some assumptions of subellipticity for the principal symbol of the operator both locally near critical points, say here $X_0=0$,
\begin{equation}\label{kps4}\inc
\exists \eps_0>0, \ \textrm{Re }p_0(X)+\eps_0H_{\textrm{Im}p_0}^2\textrm{Re }p_0(X) \sim |X|^2, \num
\end{equation}
and at infinity. In the present work, we shall not consider such a general situation where ellipticity may fail both locally and at infinity. Indeed, the main purpose of the present work being to weaken the assumptions of subellipticity near critical points, we shall simplify parts of the analysis led in \cite{HeSjSt} by requiring a property of ellipticity at infinity for the real part of the principal symbol $p_0$, but we shall consider weaker local assumptions on the doubly characteristic set. Indeed, our assumption of partial ellipticity along the singular spaces for the quadratic approximations of the principal symbol at doubly characteristic points weakens the subelliptic assumption (\ref{kps4}) since, as we shall see in Section~\ref{quadratic}, this subelliptic assumption (\ref{kps4})
induces that the singular space $S$ associated to the Hessian $q$ of the principal symbol $p_0$ at $X_0=0$ is equal to $\{0\}$. More precisely, one can check that the assumption (\ref{kps4}) is actually equivalent to the fact that the singular space $S$ is equal to zero after the intersection of exactly two kernels
\begin{equation}\label{kps5}\inc
S=\textrm{Ker}(\textrm{Re }F) \cap \textrm{Ker}\big[\textrm{Re }F(\textrm{Im }F)\big] \cap \rr^4=\{0\}.\num
\end{equation}
We refer the reader to \cite{karel} for a complete discussion of subelliptic properties of quadratic differential operators where this link between conditions (\ref{kps4}) and (\ref{kps5}) is explained.
Let us finally end this paragraph by mentioning that if one is interested in establishing semiclassical resolvent estimates for the operator $P$ instead of semiclassical hypoelliptic a priori estimates as the ones proved in Theorem~\ref{theo}, one can actually deduce resolvent estimates from that type of a priori estimates in some specific cases. This is discussed by F.~H\'erau, J.~Sj\"ostrand and C.~Stolk in~\cite{HeSjSt}, and we naturally refer the reader to~\cite{HeSjSt} (Section~11.1) for more details about this topic.

\medskip

The plan of this paper is organized as follows. Section~\ref{weight} is devoted to the construction of a global bounded weight function. Following \cite{HeSjSt}, we then recall in Section~\ref{fbi} some basic facts about the FBI-Bargmann transform and weighted spaces of holomorphic functions associated to this bounded weight function. In Section~\ref{quadratic}, we investigate the properties of the differential operators obtained by the Weyl quantization of the quadratic approximations of the principal symbol at doubly characteristic points. This study will allow us to establish in Section~\ref{tiny} some local resolvent estimates in a tiny neighborhood of these doubly characteristic points. After proving other local resolvent estimates in the exterior region (Section~\ref{intermediate}), we finally complete our proof of Theorem~\ref{theo} in Section~\ref{proof}.

\bigskip

\noindent
\textit{Remark}. We are planning to investigate in a future work the precise semiclassical asymptotics of the spectrum (modulo $\mathcal{O}(h^{\infty})$ when $h \rightarrow 0^+$) of the doubly characteristic operator $P^w(x,hD_x;h)$. More specifically, we shall try to establish  under the assumptions of Theorem~\ref{theo}, a similar result as the one proved in~\cite{HeSjSt} (Theorem~1.3) for Kramers-Fokker-Planck operators.

\bigskip

\noindent
\textit{Example.} Let $V$ and $W$ be two $C^{\infty}_b(\real^2,\rr)$ functions such that the non-negative function $V \geq 0$ is elliptic at infinity
$$\exists C >0, \forall \ |x| \geq C, \ V(x)\geq \frac{1}{C},$$
and vanishes only when $x=0$. We assume that
$$V(x)=x_1^2+\mathcal{O}(x^3),$$
while
$$W(x)=\alpha x_1^2 + 2\beta x_1 x_2 +\gamma x_2^2+\mathcal{O}(x^3),$$
when $x\rightarrow 0$, for some constants $\alpha$, $\beta$, $\gamma\in \real$, not all equal to zero. Considering the principal symbol
$$p_0(x,\xi)=\xi^2+V(x)+iW(x),$$
we notice that
$$(\textrm{Re }p_0)^{-1}(0)=\{(0,0,0,0)\},$$
and that this symbol satisfies all the assumptions (\ref{eq1.4}), (\ref{eq1.5}), (\ref{eq1.6}) and (\ref{eq1.6.5}) of Theorem~\ref{theo} with $m(x,\xi)=\langle{\xi}\rangle^2$.
The quadratic approximation of the principal symbol $p_0$ at $(0,0,0,0)$ is then given by the following quadratic form
\begin{equation}\label{eq1.7.1}\inc
q(x_1,x_2,\xi_1,\xi_2)=\xi_1^2 + \xi_2^2+ x_1^2 + i(\alpha x_1^2 + 2\beta x_1 x_2 +\gamma x_2^2),\num
\end{equation}
which is globally elliptic precisely when $\gamma \neq 0$. In general, a direct computation using (\ref{10}) and (\ref{h1}) shows that the singular space $S$ associated to $q$ is reduced to zero
precisely when $\beta^2+\gamma^2\neq 0$. In particular, when $\gamma=0$ and $\beta\neq 0$, the quadratic approximation $q$ is not globally elliptic but it obviously fulfills the assumption of partial ellipticity along its singular space $S=\{0\}$. Theorem~\ref{theo} can therefore be applied to any operator $P^w(x,hD_x;h)$ whose symbol $P(x,\xi;h)$ satisfies the following semiclassical asymptotic expansion
$$P(x,\xi;h) \sim \sum_{j=0}^{+\infty}{p_j(x,\xi)h^j},$$
with $p_j \in S(m)$ for $j \geq 1$,
despite the lack of global ellipticity of the quadratic form $q$. Finally, in the case when $\beta=\gamma=0$, the singular space $S$ is then a one-dimensional subspace and the quadratic form $q$ fails ellipticity on $S$. One can actually check in this case that the quadratic form $q$ vanishes identically on its singular space and notice that the spectrum of the associated operator
$$q^w(x,D_x)=D_{x_1}^2 + D_{x_2}^2+ (1+i\alpha)x_1^2,$$
is no longer discrete.

\bigskip

We shall finish this introduction by explaining that it is actually sufficient to establish Theorem~\ref{theo} in the special case when $m=1$. Indeed, when assuming that Theorem~\ref{theo} has already been proved when $m=1$, we may consider an order function $m \geq 1$ as in (\ref{eq1.1}) such that $m \in S(m)$; and a symbol $P(x,\xi;h)$ satisfying the associated assumptions of Theorem~\ref{theo}. Then, one can choose a symbol $\widetilde{p}_0\in S(1)$ with a non-negative real part $\textrm{Re } \widetilde{p}_0\geq 0$ which is elliptic near infinity in the symbol class $S(1)$; and such that $\widetilde{p}_0=p_0$ on a large compact set containing $p_0^{-1}(0)$ where $p_0$ stands for the principal symbol of $P(x,\xi;h)$. This is for instance the case when taking $\chi_0\in C^{\infty}_0(\real^{2n};[0,1])$ such that $\chi_0=1$ near $p_0^{-1}(0)$ and setting
$$\widetilde{p}_0=\chi_0 p_0 + (1-\chi_0).$$
Defining also the symbols
$$\widetilde{p}_j=\chi_0 p_j + (1-\chi_0) \in S(1),$$
when $j \geq 1$, we may choose $\chi\in C^{\infty}_0(\real^{2n},[0,1])$ such that $\chi=1$ near $p_0^{-1}(0)$ and $\chi_0=1$ near $\rm{supp}\, \chi$. By setting $P=P^w(x,hD_x;h)$ and $\tilde{P}=\tilde{P}^w(x,hD_x;h)$, where
$$\tilde{P}(x,\xi;h) \sim \sum_{j=0}^{+\infty}{\tilde{p}_j(x,\xi)h^j},$$
in the symbol class $S(1)$;
and using $L^2$--norms throughout, we deduce from the semiclassical elliptic regularity that
\inc\begin{align*}\label{eq1.8}
h\|u\| \leq & \  h\|\chi^w(x,hD_x) u\|+h\|(1-\chi)^w(x,hD_x)u\| \num \\
\leq & \ h\|\chi^w(x,hD_x) u\|+\mathcal{O}(h) \|(P-hz)u\|+\mathcal{O}(h^{\infty})\|u\|,
\end{align*}
since the principal symbol $p_0$ of the operator $P$ is elliptic near the support of the function $1-\chi$. By using that Theorem~\ref{theo} is valid when $m=1$, we may apply it to the operator $\widetilde{P}$ to get that if $z$ is as in Theorem~\ref{theo},
\inc\begin{align*}\label{eq1.9}
h \|\chi^w(x,hD_x) u\| \leq & \ \mathcal{O}(1)\|(\widetilde{P}-hz)\chi^w(x,hD_x) u\| \num \\
\leq & \  \mathcal{O}(1)\|(P-hz)\chi^w(x,hD_x) u\|+\mathcal{O}(h^{\infty})\|u\|,
\end{align*}
since $(\widetilde{P}-P)\chi^w(x,hD_x) = \mathcal{O}(h^{\infty})$ in $\mathcal{L}(L^2)$ when $h \rightarrow 0^+$.
We get that
\begin{equation}\label{eq1.10}\inc
h \|\chi^w(x,hD_x) u\| \leq \mathcal{O}(1)\|(P-hz)u\|+\mathcal{O}(1)\|[P,\chi^w(x,hD_x)]u\|+\mathcal{O}(h^{\infty})\|u\|.\num
\end{equation}
When estimating the commutator term in the right hand side of (\ref{eq1.10}), we take $\widetilde{\chi}\in C^{\infty}_0(\real^{2n},[0,1])$ such that $\widetilde{\chi}=1$ near $p_0^{-1}(0)$ and $\chi=1$ near $\rm{supp}\, \widetilde{\chi}$. Then, by using that
$$[P,\chi^w(x,hD_x)]\widetilde{\chi}^w(x,hD_x)=\mathcal{O}(h^{\infty}),$$
in $\mathcal{L}(L^2)$, together with the fact that $p_0$ is elliptic near the support of $1-\widetilde{\chi}$, we get that
\begin{equation}\label{eq1.11}\inc
h \|\chi^w(x,hD_x) u\| \leq \mathcal{O}(1)\|(P-hz)u\| + \mathcal{O}(h^{\infty})\|u\|,\num
\end{equation}
which in view of (\ref{eq1.8}) completes the proof of the reduction to the case when $m=1$. In what follows, we shall therefore be concerned exclusively with the case when $m=1$.

\bigskip

\noindent
\textbf{Acknowledgments:} The first author is grateful to the partial support of the National Science Foundation under grant DMS-0653275 and the Alfred P. Sloan Research Fellowship.

\section{Construction of a bounded weight function}\label{weight}
\init

The purpose of this section is to achieve the construction of a bounded weight function whose properties will be summarized below in Proposition~\ref{prop1}. When assuming that the assumptions of Theorem~\ref{theo} are all fulfilled and beginning our construction of this weight function, we shall first work in a small \neigh{} of a fixed doubly characteristic point of the principal symbol, say for example $X_1 \in p_0^{-1}(0)$. We shall assume, for notational simplicity only, that $X_1=(0,0) \in \rr^{2n}$; and drop the index 1 by denoting simply $q$ the quadratic approximation of the principal  symbol $p_0$ at $(0,0)$ appearing in (\ref{eq1.6.6}) and $S$ its associated singular space. We may therefore write that
\begin{equation}\label{eq2.0}\inc
p_0(X)=q(X)+\mathcal{O}(X^3),\num
\end{equation}
when $X \rightarrow 0$; and recall that under the assumptions of Theorem~\ref{theo}, the quadratic form $q$ is assumed to have a non-negative real part, $\textrm{Re }q \geq 0$, and to be elliptic along its singular space $S$.
Under these two assumptions, we established in \cite{HiPr} (see Section~1.4.1 and Proposition~2.0.1) that the singular space $S$ of the quadratic form $q$ has necessarily a symplectic structure and that new symplectic linear coordinates
$$X=(x,\xi)=(x',x'';\xi',\xi'') \in \real^{2n}=\real^{2n'+2n''},$$
can be chosen such that $(x'',\xi'')$ and $(x',\xi')$ are, respectively, some linear symplectic coordinates in $S$ and its symplectic orthogonal space $S^{\sigma \perp}$, so that in these coordinates, the symbol $q$ can be decomposed as the sum of two quadratic forms
\begin{equation}\label{k20b}\inc
q(x,\xi)=q|_{S^{\sigma \perp}}(x',\xi')+q|_{S}(x'',\xi''), \num
\end{equation}
where the average of the real part of first one by the flow defined by the Hamilton vector field of its imaginary part
\begin{equation}\label{k21b}\inc
\langle \textrm{Re }q|_{S^{\sigma \perp}}\rangle_{T,\textrm{Im}q|_{S^{\sigma \perp}}}(X')=\frac{1}{T}\int_{0}^{T}{\textrm{Re }q|_{S^{\sigma \perp}}(e^{tH_{\textrm{Im}q|_{S^{\sigma \perp}}}}X')dt}, \num
\end{equation}
with $X'=(x',\xi') \in \real^{2n'}$, is a positive definite quadratic form for all $T>0$; and
\begin{equation}\label{k22b}\inc
q|_S(x'',\xi'')=i \tilde{\eps}_0 \sum_{j=1}^{n''}{\lambda_j(\xi_j''^2+x_j''^2)}, \num
\end{equation}
with $\tilde{\eps}_0 \in \{\pm 1\}$, $0 \leq n'' \leq n$ and $\lambda_j>0$ for all $j=1,...,n''$.
Here the notation
$$H_f=\frac{\partial f}{\partial \xi}\cdot \frac{\partial}{\partial x} - \frac{\partial f}{\partial x} \cdot \frac{\partial}{\partial \xi},$$
stands for the Hamilton vector field of a $C^1(\real^{2d},\cc)$ function $f$. More specifically, we checked in the proof of Proposition~2.0.1 in \cite{HiPr} that the two subspaces $S$ and $S^{\sigma \perp}$ are stable by the real and imaginary parts $\textrm{Re }F$ and $\textrm{Im }F$ of the Hamilton map of the symbol $q$. By using that the flow defined by the Hamilton vector field of $\textrm{Im }q$ is the linear transformation, $e^{tH_{\textrm{Im}q}}X=e^{2t\textrm{Im}F}X$, since a direct computation using (\ref{10}) shows that $H_{\textrm{Im}q}=2\textrm{Im }F$, we deduce from (\ref{k20b}), (\ref{k21b}) and (\ref{k22b}) that
\begin{equation}\label{chel3}\inc
\textrm{Re }q(e^{tH_{\textrm{Im}q}}X)= \textrm{Re }q|_{S^{\sigma \perp}}(e^{tH_{\textrm{Im}q|_{S^{\sigma \perp}}}}X')  \num
\end{equation}
and
\begin{equation}\label{chel2}\inc
\langle \textrm{Re }q\rangle_{T,\textrm{Im}q}(X)\stackrel{\textrm{def}}{=}\frac{1}{T}\int_{0}^{T}{\textrm{Re }q(e^{tH_{\textrm{Im}q}}X)dt}=
\langle \textrm{Re }q|_{S^{\sigma \perp}}\rangle_{T,\textrm{Im}q|_{S^{\sigma \perp}}}(X'), \num
\end{equation}
where $X=(X',X'')$, $X' \in S^{\sigma \perp}$, $X'' \in S$.
We shall now use the following general observation:

\medskip

\begin{prop}\label{prop0}
For each fixed $T>0$, we have
$$\langle{\emph{\textrm{Re }}p_0\rangle}_{T,\emph{\textrm{Im}} p_0}(X)\stackrel{\emph{\textrm{def}}}{=}\frac{1}{T} \int_{0}^{T} \emph{\textrm{Re }}p_0(e^{tH_{\emph{\textrm{Im}}p_0}}X)dt =
\langle{\emph{\textrm{Re }}q\rangle}_{T,\emph{\textrm{Im}}q}(X)+\mathcal{O}(X^3),$$
when $X\rightarrow 0$; where
$$\langle{\emph{\textrm{Re }}q\rangle}_{T,\emph{\textrm{Im}} q}(X)=\frac{1}{T}\int_{0}^{T}\emph{\textrm{Re }} q(e^{tH_{\emph{\textrm{Im}}q}}X)\,dt.$$
\end{prop}

\medskip

\begin{proof}
Let $T>0$ be a fixed constant. We begin by noticing that there exists $c>0$ such that for all $0 \leq t \leq T$ and $X \in \rr^{2n}$,
\begin{equation}\label{eq2.6}\inc
|e^{tH_{\textrm{Im}p_0}}X-X| \leq c t|X| \textrm{ and } |e^{tH_{\textrm{Im}q}}X-X| \leq ct|X|. \num
\end{equation}
Indeed, there exists $C>0$ such that for all $X \in \rr^{2n}$,
\begin{equation}\label{serena1}\inc
|H_{\textrm{Im}p_0}(X)| \leq C|X| \textrm{ and } |H_{\textrm{Im}q}(X)| \leq C|X|,\num
\end{equation}
since $H_{\textrm{Im}q}=2\textrm{Im }F$, $p_0 \in S(1)$ and that $X=0$ is a doubly characteristic point of $p_0$. By writing that
\begin{equation}\label{ju10}\inc
e^{tH_{\textrm{Im}p_0}}X=X+\int_0^tH_{\textrm{Im}p_0}(e^{sH_{\textrm{Im}p_0}}X)ds, \num
\end{equation}
we notice that
$$|e^{tH_{\textrm{Im}p_0}}X| \leq |e^{tH_{\textrm{Im}p_0}}X-X|+|X| \leq |X| +C\int_0^t{|e^{sH_{\textrm{Im}p_0}}X|ds},$$
induces that
\begin{equation}\label{ju11}\inc
|e^{tH_{\textrm{Im}p_0}}X| \leq e^{Ct}|X|, \ t \geq 0, \ X \in \rr^{2n},\num
\end{equation}
by Gronwall's Lemma. By coming back to (\ref{ju10}), we easily obtain from (\ref{serena1}) the first estimate in (\ref{eq2.6}), the second one being obtained using exactly the same arguments. Then, it directly follows from (\ref{eq2.0}) and (\ref{ju11}) that
$$\langle{\textrm{Re }p_0\rangle}_{T,\textrm{Im}p_0}(X)=\langle{\textrm{Re }q\rangle}_{T,\textrm{Im}p_0}(X)+\mathcal{O}(X^3),$$
so we only need to compare the two flow averages $\langle{\textrm{Re }q\rangle}_{T,\textrm{Im}p_0}$ and $\langle{\textrm{Re }q\rangle}_{T,\textrm{Im}q}$. When doing so, it is sufficient to argue at the level of formal Taylor
expansions. Starting from (\ref{eq2.0}) and writing
\begin{equation}\label{eq2.4}\inc
\textrm{Im } p_0\sim \sum_{j=1}^{+\infty} \textrm{Im }p_{0,j}, \num
\end{equation}
where for any $j \in \nn^*$, the functions $\textrm{Im }p_{0,j}$ are homogeneous of degree $j+1$ in the variables $X=(x,\xi)$, so that in particular we have $\textrm{Im }p_{0,1}=\textrm{Im }q$; we get that
\begin{align*}
\textrm{Re }q(e^{tH_{\textrm{Im}p_0}}X) = & \ \textrm{Re }q(X) + \sum_{k=1}^{+\infty} \frac{t^k}{k!}H_{\textrm{Im}p_0}^k \textrm{Re }q(X) \\
 = & \  \textrm{Re }q(X) + \sum_{k=1}^{+\infty}\sum_{j_1=1}^{+\infty}\ldots \sum_{j_k=1}^{+\infty} \frac{t^k}{k!} H_{\textrm{Im}p_{0,j_1}}\ldots H_{\textrm{Im}p_{0,j_k}} \textrm{Re }q(X).
\end{align*}
Here, we are only interested in terms that are homogeneous of degree 2. By noticing that when $f$ is a homogeneous function of degree $j$ then $H_{\textrm{Im}p_{0,\ell}}f$ is a homogeneous function of degree $\ell+j-1$, it follows that the term $H_{\textrm{Im}p_{0,j_1}}\ldots H_{\textrm{Im}p_{0,j_k}}\textrm{Re }q$ is a homogeneous function of degree $j_1+\ldots\,+j_k+2-k$, which has a degree equal to 2 precisely when $j_1+\ldots \, +j_k=k$. We conclude that the quadratic contributions in the term $\textrm{Re }q(e^{tH_{\textrm{Im}p_0}}X)$ only come from the terms when $j_1=\ldots \,=j_k=1$,
$$\frac{t^k}{k!} H_{\textrm{Im}p_{0,1}}^k \textrm{Re }q,$$
for $k\geq 0$. This proves that
$$\textrm{Re }q(e^{tH_{\textrm{Im}p_0}}X)=\textrm{Re }q(e^{tH_{\textrm{Im}q}}X)+\mathcal{O}(X^3),$$
and ends the proof of Proposition~\ref{prop0}.
\end{proof}

\bigskip

\noindent
\textit{Remark.} Alternatively, the statement in Proposition~\ref{prop0} is an easy consequence of (\ref{eq2.6}) and the observation that there exists $\tilde{c}>0$ such that for all $0 \leq t \leq T$ and $X \in \rr^{2n}$,
\begin{equation}\label{eq2.5}\inc
|e^{tH_{\textrm{Im}p_0}}X-e^{tH_{\textrm{Im}q}}X| \leq \tilde{c}t |X|^2, \num
\end{equation}
that we shall have the occasion to use directly later on. Setting $r=\textrm{Im }p_0-\textrm{Im }q$ and writing that
\begin{align*}
e^{tH_{\textrm{Im}p_0}}X-e^{tH_{\textrm{Im}q}}X= & \ \int_0^t\big[H_{\textrm{Im}p_0}(e^{sH_{\textrm{Im}p_0}}X)-H_{\textrm{Im}q}(e^{sH_{\textrm{Im}q}}X)\big]ds \\
= & \ \int_0^tH_{r}(e^{sH_{\textrm{Im}p_0}}X)ds+\int_0^t H_{\textrm{Im}q}(e^{sH_{\textrm{Im}p_0}}X-e^{sH_{\textrm{Im}q}}X)ds,
\end{align*}
since $H_{\textrm{Im}q}$ is a linear map, the estimate (\ref{eq2.5}) directly follows from another use of Gronwall's Lemma together with (\ref{serena1}), (\ref{ju11}) and the fact that $H_r(X)=\mathcal{O}(X^2)$ uniformly on $\rr^{2n}$.

\bigskip

We can therefore deduce from (\ref{chel2}) and Proposition~\ref{prop0} that for each fixed $T>0$,
\begin{equation}\label{eq2.7}\inc
\langle{\textrm{Re }p_0\rangle}_{T,\textrm{Im}p_0}(X) = \widetilde{q}(X')+\mathcal{O}(X^3), \num
\end{equation}
where $\widetilde{q}(X')\stackrel{\textrm{def}}{=}\langle{\textrm{Re }q|_{S^{\sigma \perp}}\rangle}_{T,\textrm{Im}q|_{S^{\sigma \perp}}}(X')\sim |X'|^2$ is the positive definite quadratic form defined in (\ref{k21b}).

Let us now begin our construction of the weight function in a \neigh{} of the point $(0,0)$. When doing so, we shall follow an idea of~\cite{HeHiSj} (Section~4), and consider $g$ a decreasing $C^{\infty}(\rr_+, [0,1])$ function satisfying
\begin{equation}\label{eq2.75}\inc
g(t)=1, \ t\in [0,1]; \textrm{ and } g(t)=t^{-1}, \ t\geq 2. \num
\end{equation}
Notice that this choice induces that for each $k\in \nat$,
\begin{equation}\label{eq2.76}\inc
g^{(k)}(t)=\mathcal{O}(\langle{t\rangle}^{-1-k}),\num
\end{equation}
when $t \rightarrow +\infty$, where $\langle t \rangle=(1+t^2)^{1/2}$.
Setting
\begin{equation}\label{eq2.8}\inc
(\textrm{Re }p_0)_{\eps}(X)=g\Big(\frac{|X|^2}{\eps}\Big) \textrm{Re }p_0(X),\num
\end{equation}
for any $\eps >0$, and recalling that $p_0 \in S(1)$; we easily see from (\ref{eq2.0}) that we have the following uniform bound on $\rr^{2n}$,
\begin{equation}\label{eq2.80}\inc
(\textrm{Re }p_0)_{\eps}(X)=\mathcal{O}(\eps),\num
\end{equation}
when $\eps \rightarrow 0$.
Recalling now the well-known inequality
\begin{equation}\label{eq2.81}\inc
|f'(x)|^2 \leq 2f(x)\|f''\|_{L^{\infty}(\rr)},\num
\end{equation}
fulfilled by any non-negative smooth function with bounded second derivative, we notice that the estimate  $|\nabla \textrm{Re }p_0| =\mathcal{O}\big((\textrm{Re }p_0)^{1/2}\big)$ induces that
\begin{equation}\label{eq2.9}\inc
\partial_{X}^{\alpha}(\textrm{Re }p_0)_{\eps}(X)=\mathcal{O}(\eps^{1/2}),\num
\end{equation}
for any $\alpha\in \nat^{2n}$ with $|\alpha|=1$. By noticing from (\ref{eq2.75}) that
\begin{equation}\label{eq2.9b}\inc
(\textrm{Re }p_0)_{\eps}=\textrm{Re } p_0, \num
\end{equation}
in the region where $|X|\leq \eps^{1/2}$,
we find from (\ref{eq2.0}) that the bound (\ref{eq2.9}) improves to
\begin{equation}\label{eq2.10}\inc
\forall \alpha \in \nat^{2n}, \  |\alpha|=1, \  \partial_{X}^{\alpha} (\textrm{Re }p_0)_{\eps}(X)=\mathcal{O}(X),\num
\end{equation}
there. One can then check furthermore that for all $\alpha\in \nat^{2n}$, with $\abs{\alpha}=2$, we have
\begin{equation}\label{eq2.11}\inc
\partial_X^{\alpha} (\textrm{Re }p_0)_{\eps}(X)=\mathcal{O}(1), \num
\end{equation}
uniformly on $\rr^{2n}$.

\bigskip

\noindent
\textit{Remark}. Bounds on higher derivatives can easily be derived. In the region where $|X|\leq \eps^{1/2}$, we have
\begin{equation}\label{eq2.12}\inc
\partial_X^{\alpha} (\textrm{Re }p_0)_{\eps}(X)=\mathcal{O}\big(X^{(2-|\alpha|)_+}\big), \num
\end{equation}
for any  $\alpha\in \nat^{2n}$; while in the region where $|X| \geq \eps^{1/2}$, we check that
\begin{equation}\label{eq2.13}\inc
\partial_X^{\alpha} (\textrm{Re }p_0)_{\eps}(X)=\mathcal{O}(\eps^{1-|\alpha|/2}).\num
\end{equation}

\bigskip

\noindent
For $T>0$, we define
\begin{equation}\label{eq2.14}\inc
G_{\eps}(X) = -\int_{-\infty}^{+\infty} J\Big(-\frac{t}{T}\Big)   (\textrm{Re }p_0)_{\eps}(e^{tH_{\textrm{Im}p_0}}X)dt, \num
\end{equation}
where $J$ stands for a compactly supported piecewise affine function solving the equation
$$J'(t)=\delta(t)-\un_{[-1,0]}(t)$$
and $\un_{[-1,0]}$ the characteristic function of the set $[-1,0]$.
A direct computation as in~\cite{HeHiSj} (Section~4) using an integration by parts gives that
\begin{equation}\label{eq2.15}\inc
H_{\textrm{Im}p_0} G_{\eps} =\langle{ (\textrm{Re }p_0)_{\eps} \rangle}_{T,\textrm{Im}p_0}-(\textrm{Re }p_0)_{\eps},\num
\end{equation}
where
$$\langle (\textrm{Re }p_0)_{\eps} \rangle_{T,\textrm{Im}p_0}(X)=\frac{1}{T}\int_{0}^{T}{(\textrm{Re }p_0)_{\eps}(e^{tH_{\textrm{Im}p_0}}}  X)dt.$$
We notice from (\ref{eq2.80}) and (\ref{eq2.14}) that $G_{\eps}=\mathcal{O}(\eps)$; and that the estimates (\ref{eq2.9}) and (\ref{eq2.11}) hold for the derivatives of the function $G_{\eps}$ as well. Let us also notice from (\ref{eq2.6}), (\ref{eq2.75}), (\ref{eq2.8}) and (\ref{eq2.14}) that we have in the region, where $|X|^2 \leq \eps/2$, that
\begin{equation}\label{eq2.15.1}\inc
G_{\eps}(X) = -\int_{-\infty}^{+\infty}J\Big(-\frac{t}{T}\Big) \textrm{Re }p_0(e^{tH_{\textrm{Im}p_0}}X)dt, \num
\end{equation}
fulfills
\begin{equation}\label{chel5}\inc
H_{\textrm{Im}p_0}G_{\eps}(X)=\langle{\textrm{Re }p_0 \rangle}_{T,\textrm{Im}p_0}(X)-\textrm{Re }p_0(X);\num
\end{equation}
while we have, where $|X|^2 \geq 8\eps$, that
\begin{equation}\label{eq2.15.2}\inc
G_{\eps}(X) = - \int_{-\infty}^{+\infty}J\Big(-\frac{t}{T}\Big) \frac{\eps \textrm{Re }p_0(e^{tH_{\textrm{Im}p_0}}X)}{|e^{tH_{\textrm{Im}p_0}}X|^2}dt, \num
\end{equation}
provided that the constant $T>0$ is chosen fixed sufficiently small.

\bigskip

\noindent
\textit{Remark.} We notice from (\ref{chel3}), (\ref{eq2.15.1}) and the proof of Proposition~\ref{prop0} that in the region, where $|X|^2\leq \eps/2$, we have
\begin{equation}\label{keq}\inc
G_{\eps}(X)=G^o(X')+\mathcal{O}(X^3),\num
\end{equation}
where
\begin{equation}\label{eq2.15.2.1}\inc
G^o(X')=-\int_{-\infty}^{+\infty}J\Big(-\frac{t}{T}\Big) \textrm{Re }q|_{S^{\sigma \perp}}(e^{tH_{\textrm{Im}q|_{S^{\sigma \perp}}}}X')dt, \num
\end{equation}
is a quadratic form in the variables $X'=(x',\xi')$.

\bigskip

When considering a constant $0<\delta \ll 1$ whose value will be chosen later on, we get that
\begin{equation}\label{eq2.15.3}\inc
\widetilde{p}_0\big(X+i\delta H_{G_{\eps}}(X)\big)=p_0(X)+i \delta H_{G_{\eps}} p_0 + \mathcal{O}\left(\delta^2 |\nabla G_{\eps}|^2\right), \num
\end{equation}
when $X\in \real^{2n}$ varies in a small \neigh{} of $0$; if we denote by $\widetilde{p}_0$ an almost analytic extension of the symbol $p_0$, which is bounded together with all its derivatives in a fixed tubular \neigh{} of $\real^{2n}$. Writing
\begin{equation}\label{eq2.16}\inc
\textrm{Re}\left(\widetilde{p}_0\big(X+i\delta H_{G_{\eps}}(X)\big)\right)=\textrm{Re }p_0+\delta H_{\textrm{Im}p_0}G_{\eps}+\mathcal{O}\left(\delta^2 |\nabla G_{\eps}|^2\right) \num
\end{equation}
and
\begin{equation}\label{eq2.17}\inc
\textrm{Im}\left(\widetilde{p}_0\big(X+i\delta H_{G_{\eps}}(X)\big)\right) = \textrm{Im }p_0+\delta H_{G_{\eps}} \textrm{Re }p_0+\mathcal{O}\left(\delta^2 |\nabla G_{\eps}|^2\right), \num
\end{equation}
we shall first consider the region, where $|X|^2\leq \eps/2$. Then, by using (\ref{eq2.7}), (\ref{chel5}), as well as the fact that from (\ref{keq}), the estimate
\begin{equation}\label{chel6}\inc
\nabla G_{\eps}(X)=\mathcal{O}(X), \num
\end{equation}
is fulfilled in this region;
we get from (\ref{eq2.16}) that
\inc\begin{align*}\label{eq2.18}
\textrm{Re}\left(\widetilde{p}_0\big(X+i\delta H_{G_{\eps}}(X)\big)\right)  = & \ \textrm{Re }p_0+\delta\big(\langle{\textrm{Re }p_0}\rangle_{T,\textrm{Im}p_0}-\textrm{Re }p_0\big)+\mathcal{O}\left(\delta^2 X^2\right)  \num \\
 = & \ (1-\delta)\textrm{Re }p_0+\delta \widetilde{q}(X')+\mathcal{O}\left(\delta |X|^3+\delta^2 |X|^2\right).
\end{align*}
By noticing that (\ref{chel6}) implies that $H_{G_{\eps}} \textrm{Re }p_0(X)=\mathcal{O}(X^2)$, since $0$ is a doubly characteristic point for the symbol $p_0$; we get from (\ref{eq2.0}), (\ref{k20b}), (\ref{eq2.17}) and (\ref{chel6}) that
$$\delta \textrm{Im}\left(\widetilde{p}_0\big(X+i\delta H_{G_{\eps}}(X)\big)\right)=\delta \textrm{Im }p_0 + \mathcal{O}\left(\delta^2 X^2\right)=\delta \textrm{Im }q + \mathcal{O}\left(\delta^2 |X|^2+\delta |X|^3\right),$$
then
\inc\begin{multline*}\label{eq2.19}
i\delta \textrm{Im}\left(\widetilde{p}_0\big(X+i\delta H_{G_{\eps}}(X)\big)\right)=i\delta \textrm{Im }q|_{S^{\sigma \perp}}(X')+i\delta \textrm{Im }q|_S(X'')\\ +\mathcal{O}\left(\delta^2 |X|^2+\delta |X|^3\right),\num
\end{multline*}
when $|X|^2\leq \eps/2$. By combining (\ref{eq2.18}) and (\ref{eq2.19}), we get that
\inc\begin{multline*}\label{eq2.20}
 \textrm{Re}\left(\widetilde{p}_0\big(X+i\delta H_{G_{\eps}}(X)\big)\right)+i\delta \textrm{Im}\left(\widetilde{p}_0\big(X+i\delta H_{G_{\eps}}(X)\big)\right) \num
 = (1-\delta)\textrm{Re }p_0\\ +\delta \big(\widetilde{q}(X')+i \textrm{Im }q|_{S^{\sigma \perp}}(X')+i \textrm{Im }q|_S(X'')\big) +\mathcal{O}\left(\delta^2 |X|^2+\delta |X|^3\right).
\end{multline*}
Let us notice that the quadratic form
$$X=(X',X'')\mapsto \widetilde{q}(X')+i \textrm{Im }q|_{S^{\sigma \perp}}(X')+i \textrm{Im }q|_S(X''),$$
is elliptic on $\real^{2n}$. This comes from the facts that on one hand, the quadratic form $\widetilde{q}$ defined in (\ref{eq2.7}) is positive definite in the variables $X'$; and that on the other hand the quadratic form $\textrm{Im }q|_{S}$ defined in (\ref{k22b}) is also obviously elliptic in the variables $X''$. Combining this with the fact that, $(1-\delta)\textrm{Re }p_0\geq 0$, if $0 <\delta \leq 1$; we obtain from (\ref{eq2.20}) that there exists a positive constant $\widehat{C}>0$ such that we have in the region
where $\abs{X}^2\leq \eps/2$,
\begin{equation}\label{eq2.21}\inc
\big|\widetilde{p}_0\big(X+i\delta H_{G_{\eps}}(X)\big)\big|\geq \frac{\delta |X|^2}{\widehat{C}}-\mathcal{O}(\delta^2 X^2)-\mathcal{O}(\delta X^3)\geq \frac{\delta \abs{X}^2}{2\widehat{C}},\num
\end{equation}
when $0<\delta \leq \delta_0$ and $0<\eps \leq \eps_0$, if the positive constants $\delta_0$ and $\eps_0$ are chosen sufficiently small. This comes from the fact that one can estimate from below the quantity
\begin{align*}
& \ \big|(1-\delta)\textrm{Re }p_0 +\delta \big(\widetilde{q}(X')+i \textrm{Im }q|_{S^{\sigma \perp}}(X')+i \textrm{Im }q|_S(X'')\big)\big|\\
=& \ \Big[\big((1-\delta)\textrm{Re }p_0 +\delta \widetilde{q}(X')\big)^2+ \delta^2 \big(\textrm{Im }q|_{S^{\sigma \perp}}(X')+\textrm{Im }q|_S(X'')\big)\big)^2\Big]^{\frac{1}{2}}
\end{align*}
by the quantity
\begin{align*}
& \ \big[\delta^2 \widetilde{q}(X')^2+ \delta^2 \big(\textrm{Im }q|_{S^{\sigma \perp}}(X')+\textrm{Im }q|_S(X'')\big)^2\big]^{\frac{1}{2}}\\
=& \ \delta |\widetilde{q}(X')+i \textrm{Im }q|_{S^{\sigma \perp}}(X')+i \textrm{Im }q|_S(X'')|,
\end{align*}
since $\widetilde{q} \geq 0$.

We now consider a region of the phase space where $X$ belongs to a fixed neighborhood of $0$ and where $\abs{X}^2 \geq 8\eps$. It follows from (\ref{eq2.0}), (\ref{eq2.15}), (\ref{eq2.16}) and (\ref{eq2.17}) that
\begin{equation}\label{eq2.22}\inc
\textrm{Re}\left(\widetilde{p}_0\big(X+i\delta H_{G_{\eps}}(X)\big)\right)=\textrm{Re }p_0+\delta \big(\langle{(\textrm{Re }p_0)_{\eps}\rangle}_{T,\textrm{Im}p_0} -(\textrm{Re }p_0)_{\eps}\big)+\mathcal{O}(\delta^2 \eps)\num
\end{equation}
and
\inc\begin{align*}\label{eq2.23}
 \textrm{Im}\left(\widetilde{p}_0\big(X+i\delta H_{G_{\eps}}(X)\big)\right)= & \ \textrm{Im }p_0+\delta H_{G_{\eps}}\textrm{Re }p_0+\mathcal{O}(\delta^2 \eps) \num \\
= & \ \textrm{Im }p_0+\mathcal{O}(\delta \eps^{1/2} |X| +\delta^2 \eps)\\
= & \ \textrm{Im }q+\mathcal{O}(\delta \eps^{1/2}|X|+|X|^3),
\end{align*}
since, according to our construction of the weight function, $\nabla G_{\eps} =\mathcal{O}(\eps^{1/2})$; while we can use that $\nabla \textrm{Re }p_0(X)=\mathcal{O}(X)$, because $0$ is a doubly characteristic point of the symbol $p_0$ and that $p_0 \in S(1)$. To understand the right hand side of (\ref{eq2.22}), we first notice from (\ref{eq2.6}), (\ref{eq2.75}), (\ref{eq2.8}) and the proof of Proposition~\ref{prop0} that
\begin{align*}
\langle{(\textrm{Re }p_0)_{\eps}\rangle}_{T,\textrm{Im}p_0}(X)=& \ \eps\frac{1}{T}\int_0^T \frac{\textrm{Re }p_0(e^{tH_{\textrm{Im}p_0}}X)}{|e^{tH_{\textrm{Im}p_0}}X|^2}dt \\
= & \ \eps \frac{1}{T}\int_0^T \frac{\textrm{Re }q(e^{tH_{\textrm{Im}q}}X)} {|e^{tH_{\textrm{Im}p_0}}X|^2}dt+\mathcal{O}(\eps X),
\end{align*}
provided that $T>0$ is chosen sufficiently small. When simplifying the denominator in the right hand side of the previous equation, we use (\ref{eq2.6}) and (\ref{eq2.5}) to get that
$$\frac{1}{|e^{tH_{\textrm{Im}p_0}}X|^2}-\frac{1}{|e^{tH_{\textrm{Im}q}}X|^2}=\mathcal{O}\Big(\frac{t\abs{X}^3}{\abs{X}^4}\Big),$$
and then obtain that
\begin{equation}\label{eq2.23.1}\inc
\langle{(\textrm{Re }p_0)_{\eps}\rangle}_{T,\textrm{Im}p_0}(X)=\eps \frac{1}{T}\int_0^T \frac{\textrm{Re }q(e^{tH_{\textrm{Im}q}}X)} {|e^{tH_{\textrm{Im}q}}X|^2}dt+\mathcal{O}(\eps X), \num
\end{equation}
since $\textrm{Re }q(e^{tH_{\textrm{Im}q}}X)$ is a quadratic form.
Considering the following non-negative homogeneous function of degree $0$,
$$f(X)=\frac{1}{T}\int_0^T \frac{\textrm{Re }q(e^{tH_{\textrm{Im}q}}X)}{|e^{tH_{\textrm{Im}q}}X|^2}dt,$$
we notice from (\ref{chel3}) and (\ref{eq2.7}) that, if $f(X)=0$ then we necessarily have $\widetilde{q}(X')=0$, and therefore $X'=0$, since $\widetilde{q}$ is a positive definite quadratic form.
We may even be more precise and observe that it follows from (\ref{eq2.6}) that
\begin{equation}\label{eq2.23.2}\inc
f(X)\geq \frac{1}{\mathcal{O}(1)}\frac{|X'|^2}{|X|^2},\ X'=(x',\xi')\in S^{\sigma \perp}. \num
\end{equation}
It follows from (\ref{eq2.22}) and (\ref{eq2.23.1}) that
\begin{equation}\label{eq2.24}\inc
\textrm{Re}\left(\widetilde{p}_0\big(X+i\delta H_{G_{\eps}}(X)\big)\right)=\textrm{Re }p_0-\delta (\textrm{Re }p_0)_{\eps}+\delta \eps f(X)+\mathcal{O}(\delta^2 \eps+\delta \eps |X|). \num
\end{equation}
On the other hand, we deduce from (\ref{k20b}) and (\ref{eq2.23}) that
\begin{align*}
 \textrm{Im}\left(\widetilde{p}_0\big(X+i\delta H_{G_{\eps}}(X)\big)\right)= & \ \textrm{Im }q(X)+\mathcal{O}\big(\delta \eps^{1/2} |X| +|X|^3\big) \\
= & \  \textrm{Im }q|_{S^{\sigma \perp}}(X')+\textrm{Im }q|_S(X'')+\mathcal{O}\big(\delta \eps^{1/2} |X| +|X|^3\big),
\end{align*}
which implies, when $\abs{X}^2 \geq 8\eps$, that
\inc\begin{multline*}\label{eq2.25}
\frac{\delta \eps}{|X|^2}\textrm{Im}\left(\widetilde{p}_0\big(X+i\delta H_{G_{\eps}}(X)\big)\right)= \frac{\delta \eps}{|X|^2}
\big(\textrm{Im }q|_{S^{\sigma \perp}}(X')+\textrm{Im }q|_S(X'')\big) \\ +\mathcal{O}(\delta^2 \eps+\delta\eps |X|). \num
\end{multline*}
It follows from (\ref{eq2.24}) and (\ref{eq2.25}) that
\inc\begin{multline*}\label{eq2.25.1}
 \textrm{Re}\left(\widetilde{p}_0\big(X+i\delta H_{G_{\eps}}(X)\big)\right)+i\frac{\delta \eps}{|X|^2}\textrm{Im}\left(\widetilde{p}_0\big(X+i\delta H_{G_{\eps}}(X)\big)\right) \num
=  \textrm{Re }p_0\\ -\delta (\textrm{Re }p_0)_{\eps}  +\frac{\delta \eps}{|X|^2}\left(|X|^2 f(X)+i \textrm{Im }q|_{S^{\sigma\perp}}(X')+i \textrm{Im }q|_S(X'')\right)+\mathcal{O}\left(\delta^2 \eps+\delta \eps|X|\right),
\end{multline*}
where the continuous function of $X\neq 0$,
$$X \mapsto |X|^2 f(X)+i \textrm{Im }q|_{S^{\sigma\perp}}(X')+i \textrm{Im }q |_S(X''),$$
is homogeneous of degree 2 and does not vanish when $|X|=1$. This comes from (\ref{k22b}) and the fact that $f(X)=0$ implies that $X'=0$. By using from (\ref{eq2.75}) and (\ref{eq2.8}) that
$$\textrm{Re }p_0-\delta (\textrm{Re }p_0)_{\eps}\geq 0,$$
when $0<\delta \leq 1$; and by possibly considering a smaller constant $0<\delta_0 \leq 1$,
we deduce from (\ref{eq2.25.1}) that there exist some positive constants $C>1$ and $\tilde{C}>1$ such that for all $0<\delta \leq \delta_0$, $0<\eps \leq \eps_0$ and $8 \eps \leq |X|^2\leq 1/C$; we have
\begin{equation}\label{eq2.26}\inc
\big|\widetilde{p}_0\big(X+i\delta H_{G_{\eps}}(X)\big)\big| \geq  \frac{\delta \eps}{\tilde{C}}.\num
\end{equation}
This follows from the fact that one can estimate from below the quantity
\begin{align*}
& \ \Big|\textrm{Re }p_0-\delta (\textrm{Re }p_0)_{\eps}+\frac{\delta \eps}{|X|^2}\big(|X|^2 f(X)+i \textrm{Im }q|_{S^{\sigma\perp}}(X')+i \textrm{Im }q|_S(X'')\big)\Big|\\
= & \ \Big[\big(\textrm{Re }p_0-\delta (\textrm{Re }p_0)_{\eps}+\delta \eps f(X)\big)^2+ \frac{\delta^2 \eps^2}{|X|^4}\big(\textrm{Im }q|_{S^{\sigma\perp}}(X')+ \textrm{Im }q|_S(X'')\big)^2\Big]^{\frac{1}{2}},
\end{align*}
by
\begin{align*}
 & \ \Big[\delta^2 \eps^2 f(X)^2+ \frac{\delta^2 \eps^2}{|X|^4}\big(\textrm{Im }q|_{S^{\sigma\perp}}(X')+ \textrm{Im }q|_S(X'')\big)^2\Big]^{\frac{1}{2}}\\
=& \ \Big|\frac{\delta \eps}{|X|^2}\big(|X|^2 f(X)+i \textrm{Im }q|_{S^{\sigma\perp}}(X')+i \textrm{Im }q|_S(X'')\big)\Big|,
\end{align*}
since $f(X) \geq 0$.

Recalling (\ref{k20b}) and (\ref{k22b}), we notice that $\tilde{\eps}_0\textrm{Im }q|_S$ is a positive definite quadratic form
\begin{equation}\label{fa1}\inc
\tilde{\eps}_0 \textrm{Im }q|_S(X'')\gtrsim |X''|^2,\num
\end{equation}
while
\begin{equation}\label{fa2}\inc
\tilde{\eps}_0\textrm{Im }q|_{S^{\sigma \perp}}(X')=\mathcal{O}(X'^2).\num
\end{equation}
Continuing to work in the region where $8 \eps \leq |X|^2 \leq 1/C$ and writing
\inc\begin{multline*}\label{eq2.24.2}
\textrm{Re}\Big(\Big(1-i\tilde{\eps}_0\frac{\delta \eta \eps}{|X|^2}\Big)\widetilde{p}_0\big(X+i\delta H_{G_{\eps}}(X)\big)\Big) =  \textrm{Re}\left(\widetilde{p}_0\big(X+i\delta H_{G_{\eps}}(X)\big)\right)
\\ +\tilde{\eps}_0\frac{\delta \eta \eps}{|X|^2} \textrm{Im}\left(\widetilde{p}_0\big(X+i\delta H_{G_{\eps}}(X)\big)\right), \num
\end{multline*}
where  $0< \eta \leq 1$ stands for a constant whose value will be chosen later on, it follows from (\ref{eq2.24}) and (\ref{eq2.25}) that
\inc\begin{multline*}\label{eq2.25.1bis}
\textrm{Re}\Big(\Big(1-i\tilde{\eps}_0\frac{\delta \eta \eps}{|X|^2}\Big)\widetilde{p}_0\big(X+i\delta H_{G_{\eps}}(X)\big)\Big) =\textrm{Re }p_0(X)-\delta (\textrm{Re }p)_{\eps}(X)\\ +\frac{\delta \eps}{|X|^2}\left(|X|^2 f(X)+ \eta \tilde{\eps}_0\textrm{Im }q|_S(X'') + \eta \tilde{\eps}_0\textrm{Im }q|_{S^{\sigma \perp}}(X')\right)
+\mathcal{O}(\delta^2 \eps+\delta \eps |X|). \num
\end{multline*}
By noticing from (\ref{eq2.23.2}), (\ref{fa1}) and (\ref{fa2}) that the degree 2 homogeneous continuous function of the variable $X=(X',X'')\neq 0$,
$$X \mapsto |X|^2 f(X) + \eta \tilde{\eps}_0\textrm{Im }q |_S(X'')+\eta \tilde{\eps}_0\textrm{Im }q|_{S^{\sigma \perp}}(X'),$$
can be estimated from below by
$$|X|^2 f(X) + \eta \tilde{\eps}_0\textrm{Im }q |_S(X'')+\eta \tilde{\eps}_0\textrm{Im }q|_{S^{\sigma \perp}}(X') \geq \frac{|X|^2}{\mathcal{O}(1)},$$
provided that the constant $0 <\eta \leq 1$ is chosen sufficiently small. By using again that
$$\textrm{Re }p_0-\delta (\textrm{Re }p_0)_{\eps}\geq 0,$$
we get from (\ref{eq2.25.1bis}) after possibly decreasing the value of the constant $0<\delta_0 \leq 1$ and increasing the value of the constant $C>1$ that for all $0<\delta \leq \delta_0$ and $0<\eps\leq \eps_0$,
\begin{equation}\label{eq2.26bis}\inc
\textrm{Re}\Big(\Big(1-i\tilde{\eps}_0\frac{\delta \eta \eps}{|X|^2}\Big)\widetilde{p}_0\big(X+i\delta H_{G_{\eps}} (X)\big)\Big) \geq \frac{\delta \eps}{\mathcal{O}(1)}, \num
\end{equation}
when $8 \eps \leq |X|^2\leq 1/C$.

Having established the estimates (\ref{eq2.21}), (\ref{eq2.26}) and (\ref{eq2.26bis}), it finally remains to consider the intermediate region where $\eps/2 \leq |X|^2 \leq 8\eps$. We notice that in this region, the estimate (\ref{eq2.22}) also holds true, while we may write by using (\ref{chel3}) and Proposition~\ref{prop0} that
\begin{align*}
\langle{(\textrm{Re }p_0)_{\eps}\rangle}_{T,\textrm{Im}p_0}(X) = & \ \frac{1}{T}\int_0^T g\Big(\frac{|e^{tH_{\textrm{Im}p_0}}X|^2}{\eps}\Big) \textrm{Re }q\left(e^{tH_{\textrm{Im}q}}X\right)dt+\mathcal{O}(X^3)\\
= & \ \frac{1}{T}\int_0^T g\Big(\frac{|e^{tH_{\textrm{Im}p_0}}X|^2}{\eps}\Big) \textrm{Re }q|_{S^{\sigma \perp}}(e^{tH_{\textrm{Im}q|_{S^{\sigma \perp}}}}X')dt+\mathcal{O}(X^3),
\end{align*}
where we have here
$$g\Big(\frac{|e^{tH_{\textrm{Im}p_0}}X|^2}{\eps}\Big) \sim 1,$$
uniformly with respect to the parameters $0<\eps \leq \eps_0$ and $t\in [0,T]$; if the positive constant $T$ is chosen sufficiently small. While using the same arguments as before, it is then easy to check that the estimates (\ref{eq2.26}) and (\ref{eq2.26bis}) also hold in this region.

To summarize our discussion so far, we have shown that there exist some positive constants $C>1$, $\tilde{C}>1$, $\tilde{c} \neq 0$, $0<\eps_0 \leq 1$ and $0<\delta_0 \leq 1$ such that the local  $C^{\infty}$ weight function
$$G_{1,\eps}\stackrel{\textrm{def}}{=}G_{\eps},$$
defined in a neighborhood of $0$
satisfies for all $0 < \eps \leq \eps_0$ and $0<\delta \leq \delta_0$ that
\begin{equation}\label{eq2.27}\inc
\abs{\widetilde{p}_0\big(X+i\delta H_{G_{1,\eps}}(X)\big)}\geq \frac{\delta}{\tilde{C}} {\rm min}(|X|^2,\eps), \num
\end{equation}
when $|X| \leq 1/C$; and
\begin{equation}\label{fan10}\inc
\textrm{Re}\Big(\Big(1-i\tilde{c}\frac{\delta \eps}{|X|^2}\Big)\widetilde{p}_0\big(X+i\delta H_{G_{1,\eps}} (X)\big)\Big) \geq \frac{\delta \eps}{\tilde{C}}, \num
\end{equation}
when $\eps^{1/2} \leq |X| \leq 1/C$.

Proceeding similarly and working near each of the other doubly characteristic points $X_j\in p_0^{-1}(0)$, we obtain locally defined weight functions
$$G_{j,\eps}\in C^{\infty}({\rm neigh}(X_j,\real^{2n})),$$
for $2\leq j\leq N$, so that the natural analogues of (\ref{eq2.27}) and (\ref{fan10}) hold for $G_{j,\eps}$ in the regions where respectively $|X-X_j|\leq 1/C$ and $\eps^{1/2} \leq |X-X_j| \leq 1/C$.
To obtain a definition of the global weight function that we shall denote again $G_{\eps}$, we consider small open sets $\Omega_j \subset \real^{2n}$, $1\leq j\leq N$, with $X_j\in \Omega_j$ and
$\overline{\Omega_j}\cap \overline{\Omega_k}=\emptyset$, when $j\neq k$.
By taking some $C^{\infty}_0(\Omega_j,[0,1])$ functions $\chi_j$ such that $\chi_j=1$ in a \neigh{} of the doubly characteristic point $X_j$, we define the global $C^{\infty}_0(\real^{2n})$ weight function
$$G_{\eps} = \sum_{j=1}^N \chi_j G_{j,\eps},$$
and by restricting our attention to a fixed open set $\Omega_j$, we notice from (\ref{eq1.5}) and (\ref{eq1.6}) that the support of the functions $\nabla \chi_j$ is contained in a region where the real part of the principal symbol $p_0$ is elliptic
$$\textrm{Re }p_0(X)\geq 1/\mathcal{O}(1).$$
We have therefore proved the following result which sums up the whole discussion led in this section.

\bigskip

\begin{prop}\label{prop1}
Let $\widetilde{p}_0$ be an almost analytic extension of the principal symbol $p_0$ of the symbol $P(x,\xi;h)$ considered in Theorem~\emph{\ref{theo}} to a tubular \neigh{} of $\real^{2n}$; bounded together with all its derivatives in this neighborhood. Under the assumptions of Theorem~\emph{\ref{theo}}, one can find some constants
$$C>1, \ \tilde{C}>1, \ 0<\delta_0 \leq 1, \ 0<\eps_0 \leq 1,\ \tilde{c} \neq 0;$$
and a $C^{\infty}_0(\real^{2n},\rr)$ weight function $G_{\eps}$ depending on the parameter $0<\eps \leq \eps_0$ and supported in a neighborhood of the doubly characteristic set
$${\rm supp} \ G_{\eps} \subset \big\{X \in \rr^{2n} : {\rm dist}\big(X,p_0^{-1}(0)\big)\leq 2/C\big\},$$
such that we have uniformly for all $0< \eps \leq \eps_0$ and $0<\delta \leq \delta_0$ that
\begin{enumerate}
\item[$(i)$] $G_{\eps}=\mathcal{O}(\eps)$, $\partial^{2}G_{\eps}=\mathcal{O}(1)$ on $\rr^{2n}$
\item[$(ii)$] $\nabla G_{\eps}=\mathcal{O}\big({\rm dist}\big(X,p_0^{-1}(0)\big)\big)$ in the region where ${\rm dist}\big(X,p_0^{-1}(0)\big) \leq \eps^{1/2}$
\item[$(iii)$] $\nabla G_{\eps}=\mathcal{O}(\eps^{1/2})$ in the region where ${\rm dist}\big(X,p_0^{-1}(0)\big)\geq \eps^{1/2}$
\item[$(iv)$] We have
$$\big|\widetilde{p}_0\big(X+i\delta H_{G_{\eps}}(X)\big)\big| \geq \frac{\delta}{\tilde{C}}\emph{\textrm{min}}\big[{\rm dist}\big(X,p_0^{-1}(0)\big)^2,\eps\big],$$
in the region where ${\rm dist}\big(X,p_0^{-1}(0)\big)\leq 1/C$
\item[$(v)$] We have
$$\emph{\textrm{Re}}\Big(\Big(1-i \tilde{c}\frac{\delta \eps}{{\rm dist}\big(X,p_0^{-1}(0)\big)^2}\Big){\widetilde{p}_0\big(X+i\delta H_{G_{\eps}}(X)\big)}\Big) \geq \frac{\delta \eps}{\widetilde{C}}$$
in the region where ${\rm dist}\big(X,p_0^{-1}(0)\big) \geq \eps^{1/2}$
\end{enumerate}
\end{prop}

\section{Review of FBI tools}\label{fbi}
\init

The purpose of this section is to recall the definition of the weighted spaces of holomorphic functions associated to the weight function $G_{\eps}$ constructed in the previous section; and the action of  the operator $P$ in these spaces. The following discussion will very much follow the corresponding discussion in~\cite{HeSjSt} (Section 3)  and will therefore be somewhat brief.

Throughout this paper, we shall work with the usual semiclassical FBI-Bargmann transform
\begin{equation}\label{eq3.1}\inc
Tu (x) = \widetilde{C} h^{-3n/4} \int_{\rr^{2n}} e^{\frac{i}{h}\varphi(x,y)} u(y)dy,\ x\in \comp^n,\ \widetilde{C}>0,\num
\end{equation}
with the phase $\varphi(x,y)=\frac{i}{2}(x-y)^2$. Associated to the FBI-Bargmann transform $T$, there is the linear canonical transformation
\begin{equation}\label{eq3.2}\inc
\kappa_T: T^*\comp^n \ni (y,\eta)\mapsto (x,\xi)=(y-i\eta,\eta)\in T^*\comp^n, \num
\end{equation}
mapping the real phase space $T^*\real^n$ onto the IR-manifold
\begin{equation}\label{eq3.2.1}\inc
\Lambda_{\Phi_0}=\Big\{ \Big(x,\frac{2}{i} \frac{\partial \Phi_0}{\partial x}(x)\Big) : x \in \cc^n \Big\},\num
\end{equation}
where
\begin{equation}\label{eq3.2.1b}\inc
\Phi_0(x) =\frac{1}{2} \left(\Im x\right)^2. \num
\end{equation}
We recall that for a suitable choice of the constant $\widetilde{C}>0$ in (\ref{eq3.1}), the transformation
\begin{equation}\label{eq3.3}\inc
T: L^2(\real^n) \rightarrow H_{\Phi_0}(\comp^n) \num
\end{equation}
is unitary; where here and in what follows, when $\Phi \in C^{\infty}(\comp^n)$ is a suitable strictly plurisubharmonic weight function  close to $\Phi_0$, we shall let $H_{\Phi}(\comp^n)$ stand for the closed subspace of $L^2(\comp^n,e^{-\frac{2\Phi}{h}}L(dx))$, consisting of functions that are entire holomorphic. The integration element $L(dx)$ stands here for the Lebesgue measure on $\comp^n$.

We have an exact version of the Egorov theorem which says that
\begin{equation}\label{eq3.4}\inc
T a^w T^{-1}  = \mathfrak{a}^w, \ a\in S(1), \num
\end{equation}
where the symbol $\mathfrak{a}\in S(\Lambda_{\Phi_0},1)$ is given by the formula
$$\mathfrak{a} = a\circ \kappa_T^{-1},$$
where $\kappa_T$ is the linear canonical transformation defined in (\ref{eq3.2}).
Here the Weyl quantization of the symbol $\mathfrak{a}\in S(\Lambda_{\Phi_0},1)$, given by
$$\mathfrak{a}^w (x,hD_x) u(x)  = \frac{1}{(2\pi h)^n} \int\!\!\!\int_{\Gamma(x)} e^{\frac{i}{h}(x-y)\cdot \theta}\mathfrak{a}\Big(\frac{x+y}{2},\theta\Big) u(y)dy d\theta,$$
where $\Gamma(x)$ is the contour $\theta  = \frac{2}{i}\frac{\partial \Phi_0}{\partial x}\left(\frac{x+y}{2}\right)$, gives rise to a uniformly bounded operator on $H_{\Phi_0}(\cc^n)$. We refer to~\cite{HeSjSt} (Section~3),  as well as to~\cite{Sj95}
for a proof of this fact, based on a contour deformation argument for an almost holomorphic extension of the symbol $\mathfrak{a}$ to a tubular \neigh{} of $\Lambda_{\Phi_0}$ in $\comp^{2n}$.

Applying the remarks above to the operator $P=P^w(x,hD_x;h)$ appearing in Theorem~\ref{theo} and setting ${P}_0=T P T^{-1}$, we get that
$${P}_0  =\mathcal{O}(1): H_{\Phi_0}(\cc^n) \rightarrow H_{\Phi_0}(\cc^n),$$
if $P_0$ stands for the operator $P_0^w(x,hD_x;h)$ whose symbol has the following semiclassical asymptotic expansion
\begin{equation}\label{xi13}\inc
P_0(x,\xi;h) \sim \sum_{j=0}^{+\infty}{\mathfrak{p}_j(x,\xi) h^j},\num
\end{equation}
where $\mathfrak{p}_j=p_j\circ \kappa_T^{-1}$, for any $j \geq 0$.

Then, associated with the weight function $G_{\eps}$ constructed in Section~\ref{weight}, there is the IR-manifold
\begin{equation}\label{eq3.4.5}\inc
\Lambda_{\delta,\eps} = \left\{X+i\delta H_{G_{\eps}}(X) : X\in \real^{2n} \right \}, \num
\end{equation}
for $0 \leq \delta \leq 1$ and $0<\eps \leq \eps_0$.
Arguing as in~\cite{HeSjSt}, we see that we have
\begin{equation}\label{eq3.5}\inc
\kappa_T\left(\Lambda_{\delta,\eps}\right)=\Lambda_{\Phi_{\delta,\eps}}\stackrel{\textrm{def}}{=}\Big\{\Big(x,\frac{2}{i} \frac{\partial \Phi_{\delta,\eps}}{\partial x}(x)\Big) : x \in \comp^{2n} \Big\},\num
\end{equation}
when $0<\delta \ll 1$; where the function $\Phi_{\delta,\eps}(x)$ is the critical value with respect to the variables $(y,\eta)$ of the following functional
\begin{equation}\label{ju2}\inc
\Phi_{\delta,\eps}(x)=\textrm{v.c.}_{(y,\eta) \in \cc^n \times \rr^n}\big(-\textrm{Im }\varphi(x,y)-(\textrm{Im }y).\eta+\delta G_{\eps}(\textrm{Re }y,\eta)\big).\num
\end{equation}
Following the discussion led in \cite{HeSjSt} (Section~3.2), one can check that
$\Phi_{\delta,\eps} \in C^{\infty}(\comp^n)$ is a strictly plurisubharmonic function satisfying
\begin{equation}\label{eq3.6}\inc
\Phi_{\delta,\eps}(x)=\Phi_0(x)+\delta G_{\eps}(\Re x, -\Im x)+\mathcal{O}(\delta^2 \eps), \num
\end{equation}
uniformly on $\cc^n$, when $0< \delta \ll 1$. Since from Proposition~\ref{prop1},
$${\rm supp} \ G_{\eps} \subset \big\{X \in \rr^{2n} : {\rm dist}\big(X,p_0^{-1}(0)\big)\leq 2/C\big\},$$
we furthermore know that this strictly plurisubharmonic function $\Phi_{\delta,\eps}$ agrees with $\Phi_0$ outside a bounded set and that
\begin{equation}\label{eq3.7}\inc
\nabla\left(\Phi_{\delta,\eps}-\Phi_0\right)=\mathcal{O}(\delta \eps^{1/2}), \num
\end{equation}
uniformly on $\cc^n$;
with $\nabla^2 \Phi_{\delta,\eps}\in L^{\infty}(\comp^n)$, uniformly with respect to the parameters $\delta$ and $\eps$. Since from another use of Proposition~\ref{prop1} and (\ref{eq3.6}), we globally have
on $\cc^n$,
\begin{equation}\label{ju6}\inc
\Phi_{\delta,\eps}-\Phi_0=\mathcal{O}(\eps), \num
\end{equation}
we notice, by choosing the small parameter $\eps$ appearing in the construction of the weight $G_{\eps}$ to be equal to $\eps= Ah$,
where $h$ is the semiclassical parameter and $A \gg 1$ a large constant to be chosen in the
following; that for each fixed $A>0$, the norms in the weighted spaces $H_{\Phi_0}(\cc^n)$ and $H_{\Phi_{\delta,\eps}}(\cc^n)$ are uniformly equivalent in the semiclassical limit $h\rightarrow 0^+$. Carrying out an additional contour deformation, as in~\cite{HeSjSt}, and using (\ref{eq3.7}), we get a bounded operator
$${P}_0 = \mathcal{O}(1): H_{\Phi_{\delta,\eps}}(\cc^n) \rightarrow H_{\Phi_{\delta,\eps}}(\cc^n),$$
given, for $u$ in $H_{\Phi_{\delta,\eps}}(\cc^n)$, by
\begin{equation}\label{eq3.8}\inc
{P}_0 u(x) = \frac{1}{(2\pi h)^n} \int\!\!\!\int_{\Gamma_{\delta,\eps}(x)} e^{\frac{i}{h}(x-y)\cdot \theta} \psi_0(x-y) {P}_0\Big(\frac{x+y}{2},\theta;h\Big) u(y)dyd\theta + Ru, \num
\end{equation}
where $\psi_0$ stands for a $C^{\infty}_0(\comp^n)$ function such that $\psi_0=1$ near $0$; and $\Gamma_{\delta,\eps}(x)$ is the contour given by
$$\theta = \frac{2}{i} \frac{\partial \Phi_{\delta,\eps}}{\partial x}\Big(\frac{x+y}{2}\Big)+it_0\overline{(x-y)},\ t_0>0.$$
We have that the operator $R$ appearing in (\ref{eq3.8}) satisfies
$$R=\mathcal{O}_A(h^{\infty}): L^2\big(\comp^n,e^{-2\Phi_{\delta,\eps}/h}L(dx)\big)\rightarrow L^2\big(\comp^n,e^{-2\Phi_{\delta,\eps}/h}L(dx)\big).$$
Let us mention that in (\ref{eq3.8}), we continue to write ${P}_0$ for an almost holomorphic extension of the symbol ${P}_0 \in S(\Lambda_{\Phi_0},1)$ in a tubular \neigh{} of $\Lambda_{\Phi_0}$, bounded together with all of its derivatives. Similarly, we shall also write $\mathfrak{p}_j$ for some almost holomorphic extensions of the symbols $\mathfrak{p}_j \in S(\Lambda_{\Phi_0},1)$ in a tubular \neigh{} of $\Lambda_{\Phi_0}$, bounded together with all of their derivatives when $j \geq 0$.

\section{The quadratic case}\label{quadratic}
\init

The main purpose of this section is to get localized resolvent estimates for quadratic operators whose Weyl symbols fulfill the same properties as the quadratic approximations~$q_j$ of the principal symbol $p_0$ of the operator considered in Theorem~\ref{theo}. We shall therefore be concerned in this section with complex-valued quadratic forms with non-negative real part
\begin{equation}\label{k1}\inc
\textrm{Re } q(X) \geq 0, \  X \in \rr^{2n}, \ n \in \nn^*,\num
\end{equation}
which enjoy a property of ellipticity on their singular spaces
\begin{equation}\label{k2}\inc
X \in S, \ q(X)=0 \Rightarrow X=0. \num
\end{equation}
The main reference about quadratic operators is the work of J.~Sj\"ostrand \cite{sjostrand} and we shall actually follow very closely the analysis relying on this basic work which is led in \cite{HeSjSt} (Section~5), to explain how in view of recent improvements in the understanding of non-elliptic quadratic operators obtained in~\cite{HiPr}, the localized resolvent estimates established in~\cite{HeSjSt} (Proposition~5.2) for quadratic operators whose symbols fulfill (\ref{k1}) and the subelliptic assumption
\begin{equation}\label{k3}\inc
\exists \eps_0,C>0, \ \Re q(X) +\eps_0 H_{\textrm{Im}q}^2\textrm{Re }q(X) \geq \frac{\eps_0}{C}|X|^2, \ X \in \real^{2n}, \num
\end{equation}
can be extended to the class of quadratic operators with symbols satisfying only (\ref{k1}) and (\ref{k2}).

Let us start by verifying that the assumption (\ref{k2}) really weakens the subelliptic assumption (\ref{k3}). This is the case because condition (\ref{k3}) actually implies that the singular space is necessarily zero.
To check this, we first notice from Lemma~2 in \cite{mz} that the Hamilton map of the quadratic form $H_{\textrm{Im}q}^2\textrm{Re }q$ is given by the double commutator
$$4[\textrm{Im }F,[\textrm{Im }F, \textrm{Re }F]].$$
One can then deduce from the definition of the singular space (\ref{h1}) and the properties of skew-symmetry of the symplectic form and Hamilton maps (\ref{12}) that the two quadratic forms $\textrm{Re }q$ and $H_{\textrm{Im}q}^2\textrm{Re }q$ identically vanish on $S$,
$$\textrm{Re }q(X)=\sigma(X, \textrm{Re }F X),$$
\begin{multline*}
H_{\textrm{Im}q}^2\textrm{Re }q(X)=\sigma\big(X,4[\textrm{Im }F,[\textrm{Im }F, \textrm{Re }F]]X\big)=-8\sigma\big((\textrm{Im }F) X,[\textrm{Im }F, \textrm{Re }F]X\big) \\
= -8 \sigma\big((\textrm{Im }F) X,(\textrm{Im }F)(\textrm{Re }F)X\big)+8\sigma\big((\textrm{Im }F) X,(\textrm{Re }F)(\textrm{Im }F)X\big),
\end{multline*}
implying that when (\ref{k3}) is fulfilled, then $S=\{0\}$.

By considering from now a complex-valued quadratic form $q$ satisfying (\ref{k1}) and (\ref{k2}), we shall study the associated quadratic operator $Q=q^w(x,hD_x)$ on the FBI-Bargmann side. By using the FBI-Bargmann transform $T$ introduced in (\ref{eq3.1}),
$$T: L^2(\real^n)\rightarrow H_{\Phi_0}(\comp^n),$$
and the Egorov property recalled in (\ref{eq3.4}), we may write
\begin{equation}\label{k7}\inc
TQu=Q_0Tu,\ u\in \mathcal{S}(\real^n), \num
\end{equation}
where $Q_0$ is the quadratic differential operator on $\comp^n$ whose Weyl symbol $\mathfrak{q}_0$ is defined by the identity
\begin{equation}\label{k8}\inc
\mathfrak{q}_0\circ \kappa_T = q,\num
\end{equation}
with $\kappa_T$ being the linear canonical transformation given in (\ref{eq3.2}).
Following~\cite{Sj95}, we recall that when realizing $Q_0$ as an unbounded operator on $H_{\Phi_0}(\comp^n)$, we may first use the contour integral representation
$$Q_0 u(x) = \frac{1}{(2\pi h)^n}\int\!\!\!\int_{\theta=\frac{2}{i}\frac{\partial \Phi_0}{\partial x}\left(\frac{x+y}{2}\right)} e^{\frac{i}{h}(x-y)\cdot \theta} \mathfrak{q}_0\Big(\frac{x+y}{2},\theta\Big)u(y)\,dy\,d\theta,$$
and then, by using that the symbol $\mathfrak{q}_0$ is holomorphic, we obtain from a contour deformation the following formula for $Q_0$ as an unbounded operator on $H_{\Phi_0}(\comp^n)$,
\begin{equation}\label{k12}\inc
Q_0 u(x) = \frac{1}{(2\pi h)^n}\int\!\!\!\int_{\theta=\frac{2}{i}\frac{\partial \Phi_0}{\partial x}\left(\frac{x+y}{2}\right)+it\overline{(x-y)}} e^{\frac{i}{h}(x-y)\cdot \theta} \mathfrak{q}_0\Big(\frac{x+y}{2},\theta\Big)u(y)\,dy\,d\theta, \num
\end{equation}
for any $t>0$. We now consider as in (\ref{eq2.15.2.1}) the real-valued quadratic weight
\begin{equation}\label{k13}\inc
G^o(X) = - \int_{-\infty}^{+\infty} J\Big(-\frac{t}{T}\Big) \textrm{Re } q(e^{t H_{\textrm{Im} q}}X)dt, \num
\end{equation}
which fulfills as in (\ref{eq2.15}) the identity
\begin{equation}\label{k14}\inc
H_{\textrm{Im}q} G^o =\langle{\textrm{Re }q}\rangle_{T,\textrm{Im}q} -\textrm{Re }q, \num
\end{equation}
where
$$\langle \textrm{Re }q \rangle_{T,\textrm{Im}q}(X)=\frac{1}{T}\int_{0}^{T} \textrm{Re } q (e^{tH_{\textrm{Im}q}}X)dt.$$
As in~\cite{HeSjSt} and~\cite{HiPr}, we shall consider an IR-deformation of the real phase space $\real^{2n}$ associated to this quadratic weight $G^o$.
Setting
\begin{equation}\label{k15}\inc
\Lambda_{\delta}=\big\{X+i\delta H_{G^o}(X): X \in \real^{2n} \big\} \subset \comp^{2n}, \num
\end{equation}
for $0 \leq \delta \leq 1$, where $H_{G^o}$ stands for the Hamilton vector field of $G^o$, we find as in~\cite{HiPr} and in the previous section that we have
\begin{equation}\label{k16}\inc
\kappa_T(\Lambda_{\delta})=\Lambda_{\Phi_{\delta}^o}\stackrel{\textrm{def}}{=}\Big\{\Big(x,\frac{2}{i}\frac{\partial \Phi_{\delta}^o}{\partial x}(x)\Big) : x\in \comp^n\Big\}, \num
\end{equation}
for all $0\leq \delta \leq \delta_0$,  with $\delta_0>0$ small enough, where $\Phi_{\delta}^o$ is a strictly plurisubharmonic quadratic form on $\comp^n$ verifying
\begin{equation}\label{k17}\inc
\Phi_{\delta}^o(x)=\Phi_0(x)+\delta G^o(\Re x, -\Im x)+\mathcal{O}(\delta^2 \abs{x}^2). \num
\end{equation}
Associated to the quadratic form $\Phi_{\delta}^o$ is the weighted space of holomorphic functions $H_{\Phi_{\delta}^o}(\comp^n)$ defined as in Section~\ref{fbi}. We can now view the operator $Q_0$ as an unbounded operator
$$Q_0: H_{\Phi_{\delta}^o}(\comp^n) \rightarrow H_{\Phi_{\delta}^o}(\comp^n),$$
if we make a new contour deformation as in (\ref{k12}) and set
\begin{equation}\label{k18}\inc
Q_0 u(x) = \frac{1}{(2\pi h)^n}\int\!\!\!\int_{\theta=\frac{2}{i}\frac{\partial \Phi_{\delta}^o}{\partial x}\left(\frac{x+y}{2}\right)+it_0\overline{(x-y)}} e^{\frac{i}{h}(x-y)\cdot \theta}\mathfrak{q}_0\Big(\frac{x+y}{2},\theta\Big)u(y)\,dy\,d\theta, \num
\end{equation}
for $t_0>0$, when $\delta$ is sufficiently small, $0 < \delta \ll 1$. By coming back to the real side via the FBI-Bargmann transform, the operator $Q_0$ can then be viewed as an unbounded operator on $L^2(\real^n)$ whose Weyl symbol is given by the following quadratic form
\begin{equation}\label{k19}\inc
\widetilde{q}(X)=q\big(X+i \delta H_{G^o}(X)\big).\num
\end{equation}

As in Section~\ref{weight}, we recall that the singular space $S$ of a quadratic symbol $q$ satisfying the assumptions (\ref{k1}) and (\ref{k2}) always has a symplectic structure (see Section~1.4.1 in~\cite{HiPr}) and that according to~\cite{HiPr} (Proposition 2.0.1), new symplectic linear coordinates
$$X=(x,\xi)=(x',x'';\xi',\xi'') \in \real^{2n}=\real^{2n'+2n''},$$
can be chosen such that $(x'',\xi'')$ and $(x',\xi')$ are, respectively, some linear symplectic coordinates in $S$ and its symplectic orthogonal space $S^{\sigma \perp}$, so that in these coordinates, the symbol $q$ can be decomposed as the sum of two quadratic forms
\begin{equation}\label{k20}\inc
q(x,\xi)=q|_{S^{\sigma \perp}}(x',\xi')+q|_{S}(x'',\xi''), \num
\end{equation}
such that the average of the real part of first one by the flow defined by the Hamilton vector field of its imaginary part
\begin{equation}\label{k21}\inc
\langle \textrm{Re }q|_{S^{\sigma \perp}}\rangle_{T,\textrm{Im}q|_{S^{\sigma \perp}}}(X')=\frac{1}{T}\int_{0}^{T}{\textrm{Re }q|_{S^{\sigma \perp}}(e^{tH_{\textrm{Im}q|_{S^{\sigma \perp}}}}X')dt},  \num
\end{equation}
where $X'=(x',\xi') \in \real^{2n'}$;
is a positive definite quadratic form for all $T>0$; and
\begin{equation}\label{k22}\inc
q|_S(x'',\xi'')=i \tilde{\eps}_0 \sum_{j=1}^{n''}{\lambda_j(\xi_j''^2+x_j''^2)}, \num
\end{equation}
with $\tilde{\eps}_0 \in \{\pm 1\}$, $0 \leq n'' \leq n$ and $\lambda_j>0$ for all $j=1,...,n''$.
Then, by noticing that a direct computation gives that
$$H_{\textrm{Im}q}=2 \textrm{Im }F,$$
we deduce from (\ref{k13}), (\ref{k14}), (\ref{k20}), (\ref{k22}); and the stability property of the spaces $S$ and $S^{\sigma \perp}$ by the map $\textrm{Im }F$ \big(see (2.0.4) in \cite{HiPr}\big); that the quadratic weight $G^o$ only depends on the variables $X'=(x',\xi') \in \real^{2n'}$,
\begin{equation}\label{k23}\inc
G^o(X')=- \int_{-\infty}^{+\infty} J\Big(-\frac{t}{T}\Big) \Re q|_{S^{\sigma \perp}}(e^{tH_{\textrm{Im} q|_{S^{\sigma \perp}}}}X')\,dt \num
\end{equation}
and satisfies
\begin{equation}\label{k24}\inc
H_{\textrm{Im} q|_{S^{\sigma \perp}}} G^o =\langle{\Re q|_{S^{\sigma \perp}}}\rangle_{T,\textrm{Im}q|_{S^{\sigma \perp}}} -\Re q|_{S^{\sigma \perp}}. \num
\end{equation}
Since the symbol $\widetilde{q}$ in (\ref{k19}) is easily seen to be equal to
\begin{equation}\label{k25}\inc
\widetilde{q}=q-i \delta H_q G^o +\mathcal{O}(\delta^2\abs{\nabla G^o}^2), \num
\end{equation}
we deduce from the previous tensorization of the variables (\ref{k20}), (\ref{k22}) and (\ref{k23}) that this symbol can be written as
\begin{equation}\label{k26}\inc
\widetilde{q}(X)=r(X')+i \tilde{\eps}_0 \sum_{j=1}^{n''}{\lambda_j(\xi_j''^2+x_j''^2)}, \num
\end{equation}
with
$$r(X')=q|_{S^{\sigma \perp}}(X')-i \delta H_{q|_{S^{\sigma \perp}}}G^o(X')+\mathcal{O}(\delta^2\abs{X'}^2).$$
Using now (\ref{k1}), (\ref{k20}), (\ref{k21}) and (\ref{k24}), we notice that the real part of the quadratic symbol $r$,
\begin{multline*}
\textrm{Re }q|_{S^{\sigma \perp}}(X')+\delta H_{\textrm{Im}q|_{S^{\sigma \perp}}}G^o(X')+\mathcal{O}(\delta^2\abs{X'}^2) \\
=(1-\delta)\textrm{Re }q|_{S^{\sigma \perp}}(X')+\delta \langle \textrm{Re }q|_{S^{\sigma \perp}}\rangle_{T,\textrm{Im}q|_{S^{\sigma \perp}}}(X') +\mathcal{O}(\delta^2\abs{X'}^2)
\geq \frac{\delta}{C}|X'|^2,
\end{multline*}
for $C>1$; is a positive definite quadratic form for all $0<\delta \ll 1$ sufficiently small. In view of (\ref{k26}), this particular  property implies the ellipticity of the quadratic symbol $\widetilde{q}$ on $\real^{2n}$. We can then apply the classical result of J.~Sj\"ostrand (Theorem 3.5 in \cite{sjostrand}) to the operator $Q_0$ viewed as an unbounded operator on $H_{\Phi_{\delta}^o}(\comp^n)$, and using similar arguments as the ones used by F.~H\'erau, J.~Sj\"ostrand and C.~Stolk in their proof of Proposition 5.1 in \cite{HeSjSt}, we get that this operator $Q_0$ fulfills all the properties stated in~\cite{HeSjSt} (Proposition~5.1), namely that its spectrum, as an unbounded operator on $H_{\Phi_{\delta}^o}(\comp^n)$ for $0<\delta \ll 1$, is only composed of eigenvalues of finite multiplicity with the following structure
\begin{equation}\label{k27}\inc
\sigma\big{(}Q_0\big{)}=\Big\{ \sum_{\substack{\lambda \in \sigma(F), \\  -i \lambda \in \comp_+ \cup (\Sigma(q|_S) \setminus \{0\})}}{h\big{(}r_{\lambda}+2 k_{\lambda}\big{)}(-i\lambda) : k_{\lambda} \in \nat}\Big\}, \num
\end{equation}
where $F$ is the Hamilton map associated to the quadratic form $q$, $r_{\lambda}$ the dimension of the space of generalized eigenvectors of $F$ in $\comp^{2n}$ belonging to the eigenvalue $\lambda \in \comp$,
$$\Sigma(q|_S)=\overline{q(S)} \textrm{ and } \comp_+=\{z \in \comp: \textrm{Re }z>0\};$$
and that if $z$ remains in a compact set of empty intersection with $\sigma(Q_0)$, then we have with $d(x)=|x|$,
\begin{equation}\label{k28}\inc
\|(h+d^2)u\| \leq \mathcal{O}(1) \|(Q_0-hz)u\|, \num
\end{equation}
and
\begin{equation}\label{k29}\inc
\|(h+d^2)^{\frac{1}{2}}u\| \leq \mathcal{O}(1) \|(h+d^2)^{-\frac{1}{2}}(Q_0-hz)u\|, \num
\end{equation}
for all holomorphic function $u$ satisfying respectively
$$(h+d^2)u \in L^2\big(\comp^n,e^{-2\Phi_{\delta}^o(x)/h}L(dx)\big) \textrm{ and }(h+d^2)^{\frac{1}{2}}u \in L^2\big(\comp^n,e^{-2\Phi_{\delta}^o(x)/h}L(dx)\big).$$
The only difference with Proposition~5.1 in \cite{HeSjSt} is that we describe here the spectrum of $Q_0$ by the mean of the singular space as we did in \cite{HiPr} (Theorem 1.2.2).

We have therefore checked that any quadratic operator whose symbol $q$ fulfills the assumptions (\ref{k1}) and (\ref{k2}), verifies all the properties stated in~\cite{HeSjSt} (Proposition~5.1). We can then deduce from Proposition~5.2 in \cite{HeSjSt} that the operator $Q_0$ also verifies the estimates stated in~\cite{HeSjSt} (Proposition~5.2), since the proof of Proposition~5.2 only relies on estimates obtained in~\cite{HeSjSt} (Proposition~5.1). This proves that any quadratic operator whose Weyl symbol fulfills the assumptions (\ref{k1}) and (\ref{k2}), defines on the FBI transform side, an unbounded operator on the spaces $H_{\Phi_{\delta}^o}(\comp^n)$, with $0 < \delta \ll 1$, which fulfills the following localized resolvent estimates:

\bigskip

\begin{lemma}\label{klem1}
Let $\chi_0 \in C_0^{\infty}(\comp^n)$ be fixed and equal to {\rm 1} near {\rm 0}, and fix $k \in \real$. Then for $z$ varying in a compact set that does not contain any eigenvalues of the operator $Q_0|_{h=1}$ described in \emph{(\ref{k27})}, we have with $d(x)=|x|$,
\begin{equation}\label{k30}\inc
\|(h+d^2)^{1-k}\chi_0 u\| \leq \mathcal{O}(1)\|(h+d^2)^{-k}\chi_0(Q_0-hz)u\| + \mathcal{O} (h^{\frac{1}{2}})\|\un_K u\|, \num
\end{equation}
where $K$ is any fixed neighborhood of $\emph{\textrm{supp}}(\nabla\chi_0)$ and $\un_K$ stands for its characteristic function. Here the norm is taken in the space $L^2(\comp^n, e^{-2\Phi_{\delta}^o/h}L(dx))$.
\end{lemma}

\bigskip

By using these localized resolvent estimates satisfied by quadratic operators defined by the quadratic approximations of the principal symbol $p_0$ at double characteristic points, we shall establish in the next section local resolvent estimates for these operators and the operator $P=P^w(x,hD_x;h)$ in a tiny neighborhood of the doubly characteristic set.

\section{Local resolvent estimates in a tiny \neigh{} of the doubly characteristic set}\label{tiny}
\init

In all of this section, $G_{\eps}$ stands for the weight function constructed in Proposition~\ref{prop1} and we choose the small parameter $\eps$ to be equal to
\begin{equation}\label{k31}\inc
\eps=Ah, \ A \gg 1, \num
\end{equation}
where $h$ is the semiclassical parameter and $A$ a large constant to be chosen in the following.

We recall from Section~\ref{fbi} that the IR-manifold $\Lambda_{\delta,\eps}$ defined in (\ref{eq3.4.5}) is represented on the FBI-transform side by
\begin{equation}\label{k32}\inc
\kappa_{T}(\Lambda_{\delta,\eps})=\Lambda_{ \Phi_{\delta,\eps}}=\Big\{\Big(x,\frac{2}{i}\frac{\partial \Phi_{\delta,\eps}}{\partial x}(x)\Big): x \in \comp^n\Big\}, \num
\end{equation}
where $ \Phi_{\delta,\eps} \in C^{\infty}(\cc^n)$ is the strictly plurisubharmonic function introduced in (\ref{ju2}) for $0<\eps \leq \eps_0$ and sufficiently small values of the parameter $0<\delta \ll 1$.
Let us notice directly from (\ref{eq2.15.1}) and (\ref{ju2}) that this weight function $\Phi_{\delta,\eps}$ is independent of the parameter $\eps$ in a region where $|x| \leq \sqrt{\eps}/C$, after a suitable choice of the constant $C>0$. By making a rescaling in $\eps$, we may and will assume in the following that we have $C=1$.

By recalling from Section~\ref{fbi} that we have a uniformly bounded operator
\begin{equation}\label{k33.1}\inc
{P}_0: H_{\Phi_{\delta,\eps}}(\cc^n) \rightarrow H_{\Phi_{\delta,\eps}}(\cc^n), \num
\end{equation}
we shall be concerned with the almost holomorphic extension of the principal symbol $\mathfrak{p}_0$ of the operator ${P}_0$, restricted to $\Lambda_{\Phi_{\delta,\eps}}$, and to ease the notation in this section, we shall simply write
$$\mathfrak{p}_0\stackrel{\textrm{def}}{=} \mathfrak{p}_0|_{\Lambda_{\Phi_{\delta,\eps}}}.$$
We shall also assume for simplicity that the characteristic set $p_0^{-1}(0)$ is composed of an unique point, say here $(0,0)$, which corresponds to the point $\comp^n\ni x=0$ on the FBI-transform side.

Considering the quadratic approximation of the principal symbol $p_0$,
\begin{equation}\label{k34}\inc
q_0(x,\xi)=\sum_{|\alpha+\beta|=2}\frac{\partial_x^{\alpha}\partial_{\xi}^{\beta}p_0(0,0)}{\alpha! \beta!}x^{\alpha}\xi^{\beta}, \num
\end{equation}
at the critical point $(0,0)$, we may write
\begin{equation}\label{k35}\inc
p_0(x,\xi)+hp_1(x,\xi)-q_0(x,\xi)-hp_1(0,0)=\mathcal{O}\big((h+|(x,\xi)|^2)^{\frac{3}{2}}\big), \ \ 0<h \leq 1, \num
\end{equation}
where $p_1$ stands for the subprincipal symbol of the symbol $P(x,\xi;h)$ appearing in its asymptotic expansion (\ref{xi1}).
Using this quadratic approximation, we aim in this section at getting local resolvent estimates for the operator $P_0$ in the following tiny neighborhood of the doubly characteristic point
$$|x| \leq \sqrt{\eps},$$
on the FBI-transform side. Working in this region, we can first notice that we may replace $\Phi_{\delta,\eps}$ by the strictly plurisubharmonic quadratic form $\Phi_{\delta}^o$, introduced in (\ref{k16}), since the identities (\ref{keq}), (\ref{eq3.6}) and (\ref{k17}) induce that the two $L^2$-norms associated to these weight functions are equivalent in the region where $|x| \leq \sqrt{\eps}$,
\begin{equation}\label{ju1}\inc
\exp\big(-\mathcal{O}(1)(A+ A^{3/2} h^{1/2})\big) \leq e^{-\frac{2\Phi_{\delta,\eps}(x)}{h}} e^{\frac{2\Phi_{\delta}^o(x)}{h}} \leq \exp\big(\mathcal{O}(1)(A+ A^{3/2}h^{1/2})\big), \num
\end{equation}
for each fixed constant $A\gg 1$.

By assuming therefore that the weight function $\Phi_{\delta,\eps}$ is equal to the quadratic form $\Phi_{\delta}^o$ in the region where $|x| \leq \sqrt{\eps}$, we realize the operators $\mathfrak{p}_0^w(x,hD_x)$, $\mathfrak{p}_1^w(x,D_x)$, $P_0$ and $Q_0$ with a contour as in (\ref{eq3.8}). We deduce from this realization and (\ref{k35}) that the difference between the corresponding effective kernels of $\mathfrak{p}_0^w(x,D_x)+h\mathfrak{p}_1^w(x,hD_x)$ and $Q_0+hp_1(0,0)$ is, for some $D>0$,
$$\mathcal{O}_A(1)h^{-n}e^{-\frac{D}{h}|x-y|^2}(h+|x|^2+|y|^2)^{\frac{3}{2}}=\mathcal{O}_A(1)h^{-n}e^{-\frac{D}{h}|x-y|^2}(h^{\frac{3}{2}}+|x|^3+|x-y|^3),$$
which implies as in~\cite{HeSjSt} (see (6.9)) that
\inc\begin{multline*}\label{k36}
\|\mathfrak{p}_0^w(x,hD_x)u+h\mathfrak{p}_1^w(x,hD_x)u-Q_0u-hp_1(0,0)u\|_{H_{\Phi_{\delta,\eps}}(|x| \leq \sqrt{Ah})}\\
=\mathcal{O}_A(h^{\frac{3}{2}})\|u\|_{H_{\Phi_{\delta,\eps}}}, \num
\end{multline*}
where the notation $\mathcal{O}_A$ is here to emphasize that the constant may depend on the large parameter $A$ which will be chosen later on.

Keeping in mind that the operator $Q_0$ is realized with a contour as in (\ref{eq3.8}) and using the fact that its symbol is a quadratic polynomial, we check that if the quantity $Q_0u$ is replaced by the corresponding differential expression
$$\sum_{|\alpha+\beta|=2}\frac{\partial_x^{\alpha}\partial_{\xi}^{\beta}p_0(0,0)}{\alpha! \beta!}\big(x^{\alpha}(hD_x)^{\beta}\big)^wu,$$
on the FBI-transform side, we then commit an error $w$ satisfying
\begin{equation}\label{k37}\inc
\|w\|_{H_{\Phi_{\delta,\eps}}(|x| \leq \sqrt{Ah})} \leq e^{-\frac{1}{Ch}} \|u\|_{H_{\Phi_{\delta,\eps}}}, \num
\end{equation}
for some $C>0$.

Continuing to follow \cite{HeSjSt} (Section~6), we can now consider the change of variables on the FBI-transform side
$$x=\sqrt{\eps}\tilde{x}, \ hD_{x}=\sqrt{\eps}\tilde{h}D_{\tilde{x}},$$
with $\tilde{h}=A^{-1}$. Then, setting $\tilde{Q}_0=\mathfrak{q}_0(\tilde{x},\tilde{h}D_{\tilde{x}})^w$, homogeneity properties imply that
\begin{equation}\label{k38}\inc
Q_0=\mathfrak{q}_0^w(x,hD_x)=\frac{h}{\tilde{h}}\mathfrak{q}_0^w(\tilde{x},\tilde{h}D_{\tilde{x}})=\frac{h}{\tilde{h}}\tilde{Q}_0, \num
\end{equation}
and, with $d=d(x)=|x|$, $\tilde{d}=d(\tilde{x})=|\tilde{x}|$; we have the following identities, since $\Phi_{\delta}^o(x)$ is a quadratic form
\begin{equation}\label{k39}\inc
h+d^2=\frac{h}{\tilde{h}}(\tilde{h}+\tilde{d}^2), \ e^{-\frac{2\Phi_{\delta}^o(x)}{h}}=e^{-\frac{2\Phi_{\delta}^o(\tilde{x})}{\tilde{h}}}. \num
\end{equation}

Let $\chi_0 \in C_0^{\infty}(\comp^n)$ be fixed and equal to 1 near the set where $\abs{x}\leq 1$. For $k \in \real$, we can then apply Lemma~\ref{klem1} to the operator $\tilde{Q}_0$ to get that for any $z$ belonging to a fixed compact set avoiding the eigenvalues of the operator $\tilde{Q}_0|_{\tilde{h}=1}$, the following estimate is fulfilled
\begin{equation}\label{k40}\inc
\|(\tilde{h}+\tilde{d}^2)^{1-k}\chi_0(\tilde{x}) \tilde{u}\| \leq  \mathcal{O}(1) \|(\tilde{h}+\tilde{d}^2)^{-k}\chi_0(\tilde{x})(\tilde{Q}_0-\tilde{h}z)\tilde{u}\| + \mathcal{O}(\tilde{h}^{\frac{1}{2}})\|\un_{K}\tilde{u}\|, \num
\end{equation}
where $K$ is a fixed neighborhood of $\textrm{supp}(\nabla \chi_0)$. Here the norm is taken with respect to the weight $e^{-2\Phi_{\delta}^o(\tilde{x})/\tilde{h}}$. By noticing that we can replace the term $\mathcal{O}(\tilde{h}^{\frac{1}{2}})\|\un_{K}\tilde{u}\|$ by
$$\mathcal{O}(\tilde{h}^{\frac{1}{2}})\|(\tilde{h}+\tilde{d}^2)^{1-k}\un_{K}\tilde{u}\|,$$
and that in view of (\ref{k38}) the two operators $Q_0|_{h=1}$ and $\tilde{Q}_0|_{\tilde{h}=1}$ have the same spectra, we get from (\ref{k38}), (\ref{k39}) and (\ref{k40}) by coming back to initial variables $x$, $h$, that for any $z$ belonging to a fixed compact set avoiding the eigenvalues of the operator $Q_0|_{h=1}$, we have
\begin{multline*}
\Big\| \Big(\frac{\tilde{h}}{h}\Big)^{1-k} (h+d^2)^{1-k}\chi_0\Big(\frac{x}{\sqrt{Ah}}\Big) u\Big\| \leq \mathcal{O}(1) \Big\|  \Big(\frac{\tilde{h}}{h}\Big)^{1-k} (h+d^2)^{-k}\chi_0\Big(\frac{x}{\sqrt{Ah}}\Big)  (Q_0-hz)u\Big\|\\ +\tilde{C}\tilde{h}^{\frac{1}{2}}   \Big\|\Big(\frac{\tilde{h}}{h}\Big)^{1-k}(h+d^2)^{1-k}\un_{K}\Big(\frac{x}{\sqrt{Ah}}\Big)u\Big\|,
\end{multline*}
which finally induces the following result:

\bigskip

\begin{lemma}\label{klem1cel}
Let $\chi_0 \in C_0^{\infty}(\comp^n)$ be fixed and equal to {\rm 1} near the set where $\abs{x}\leq 1$, and fix $k \in \real$. Then for $z$ varying in a compact set that does not contain any eigenvalues of the operator $Q_0|_{h=1}$ described in \emph{(\ref{k27})}, we have with $d(x)=|x|$,
\inc\begin{multline}\label{k41}
\Big\| (h+d^2)^{1-k}\chi_0\Big(\frac{x}{\sqrt{Ah}}\Big) u\Big\| \leq \mathcal{O}(1) \Big\| (h+d^2)^{-k}\chi_0\Big(\frac{x}{\sqrt{Ah}}\Big)  (Q_0-hz)u\Big\|\\ +\mathcal{O}\Big(\frac{1}{\sqrt{A}}\Big) \Big\|(h+d^2)^{1-k}\un_{K}\Big(\frac{x}{\sqrt{Ah}}\Big)u\Big\|,\num
\end{multline}
when $0<h \ll 1$ and $A \gg 1$,
where $K$ is any fixed neighborhood of $\emph{\textrm{supp}}(\nabla\chi_0)$ and $\un_K$ stands for its characteristic function. Here the norm is taken with respect to the weight $e^{-2\Phi_{\delta}^o(x)/h}$. According to \emph{(\ref{ju1})}, all these norms with respect to the weight  $e^{-2\Phi_{\delta}^o/h}$ can then be replaced by norms with respect to the weight $e^{-2\Phi_{\delta,\eps}/h}$ in the previous estimate.
\end{lemma}

\section{Local resolvent estimates in the exterior region}\label{intermediate}
\init

We shall now establish local resolvent estimates in the region outside the tiny neighborhood of the doubly characteristic set considered in the previous section. As in Section~\ref{tiny}, we shall continue to assume for simplicity only that the doubly characteristic set $p_0^{-1}(0)$ is composed of a unique point, say here $(0,0) \in \rr^{2n}$, which corresponds to the point $\comp^n\ni x=0$ on the FBI-transform side.

Considering always the same IR-manifold
$$\Lambda_{\Phi_{\delta,\eps}}=\kappa_{T}(\Lambda_{\delta,\eps})=\Big\{\Big(x,\frac{2}{i}\frac{\partial \Phi_{\delta,\eps}}{\partial x}(x)\Big) : x \in \comp^n\Big\},$$
associated to the weight $G_{\eps}$ as in the previous section, we recall that the small parameter $\eps$ is taken equal to
$$\eps=Ah, \ A \gg 1,$$
where $h$ stands for the semiclassical parameter and $A$ a large parameter still remaining to be chosen. The purpose of this section is to get a local resolvent estimate in the region outside the tiny $\sqrt{\eps}$--neighborhood of the doubly characteristic point $0$ studied in the previous section. We shall therefore be concerned with studying the following region on the FBI-transform side of the IR-manifold~$\Lambda_{\Phi_{\delta,\eps}}$,
\begin{equation}\label{k42}\inc
|x| \geq \sqrt{\eps}. \num
\end{equation}
We begin by noticing from (\ref{eq3.2.1b}) and (\ref{eq3.6}) that we have
\begin{equation}\label{eq6.2xi}\inc
\xi = -\textrm{Im }x +\mathcal{O}(\delta \sqrt{\eps}),\num
\end{equation}
when $(x,\xi)\in \Lambda_{\Phi_{\delta,\eps}}$.
When working in the unbounded region (\ref{k42}), we recall from Proposition~\ref{prop1} that we have
\begin{equation}\label{eq6.2}\inc
\textrm{Re}\Big(\Big(1-i \tilde{c} \frac{\delta \eps}{|x|^2}\Big) \mathfrak{p}_0\Big(x,\frac{2}{i}\frac{\partial \Phi_{\delta,\eps}(x)}{\partial x}\Big)\Big) \geq \frac{\delta \eps}{\tilde{C}},\num
\end{equation}
when $|x| \geq \sqrt{\eps}$.
It is therefore convenient to consider again the new variables
\begin{equation}\label{eq6.3}\inc
x = \sqrt{\eps} \widetilde{x}.\num
\end{equation}
In these new coordinates, the IR-manifold $\Lambda_{\Phi_{\delta,\eps}}$ then becomes replaced by
\begin{equation}\label{eq6.4}\inc
\Lambda_{\widetilde{\Phi}_{\delta,\eps}}= \Big\{\Big(\widetilde{x},\frac{2}{i} \frac{\partial \widetilde{\Phi}_{\delta,\eps}(\widetilde{x})}{\partial
\widetilde{x}},\Big) : \widetilde{x} \in \comp^n\Big\},\num
\end{equation}
with
\begin{equation}\label{eq6.4xi}\inc
\widetilde{\Phi}_{\delta,\eps}(\widetilde{x}) = \frac{1}{\eps} \Phi_{\delta,\eps}(\sqrt{\eps}\widetilde{x}).\num
\end{equation}
We notice from (\ref{eq6.2xi}) and (\ref{eq6.4xi}) that the function $\nabla^2 \widetilde{\Phi}_{\delta,\eps}$ belongs to the space $L^{\infty}(\cc^{n})$ uniformly with respect to the parameters $0<\delta \leq \delta_0$ and $0<\eps \leq \eps_0$, since it is the case for the function $\nabla^2 \Phi_{\delta,\eps}$; and that
$$\widetilde{\xi} = -\Im \widetilde{x} +\mathcal{O}(\delta),$$
along the IR-manifold $\Lambda_{\widetilde{\Phi}_{\delta,\eps}}$.
Writing
\begin{equation}\label{eq6.5}\inc
\frac{1}{\eps} \mathfrak{p}_0^w(x, hD_x) = \frac{1}{\eps}\mathfrak{p}_0^w\big(\sqrt{\eps}(\widetilde{x},\widetilde{h}D_{\widetilde{x}})\big),\num
\end{equation}
with
\begin{equation}\label{eq6.5xi}\inc
\widetilde{h}=\frac{h}{\eps}=\frac{1}{A},\num
\end{equation}
we shall work with the $\widetilde{h}$--pseudodifferential operator
$$\widetilde{P}_{\eps} = \frac{1}{\eps} \mathfrak{p}_0^w(x, hD_x),$$
whose Weyl symbol is given by
\begin{equation}\label{eq6.6}\inc
\widetilde{\mathfrak{p}}_{\eps}(\widetilde{x},\widetilde{\xi}) = \frac{1}{\eps} \mathfrak{p}_0\big(\sqrt{\eps}(\widetilde{x},\widetilde{\xi})\big), \num
\end{equation}
fulfilling for any $k\geq 0$,
$$\nabla^{k}\widetilde{\mathfrak{p}}_{\eps} = \mathcal{O}(\eps^{k/2-1}),$$
so that in particular, we have
\begin{equation}\label{eq6.6.1}\inc
\nabla^k \widetilde{\mathfrak{p}}_{\eps} = \mathcal{O}(1),\num
\end{equation}
uniformly with respect to the parameter  $0 < \eps \leq \eps_0$ on $\Lambda_{\widetilde{\Phi}_{\delta,\eps}}$, when $k\geq 2$. Recalling that
$$\mathfrak{p}_0=p_0 \circ \kappa_T^{-1},$$
where $\kappa_T$ is the linear canonical transformation defined in (\ref{eq3.2}), this remark together with the fact that we have the estimate
$$\mathfrak{p}_0(x,\xi)=\mathcal{O}(|(x,\xi)|^2),$$
near the origin since $(0,0)$ is a doubly characteristic point for the symbol $p_0$; imply that in any region along $\Lambda_{\widetilde{\Phi}_{\delta,\eps}}$ where the variables $\widetilde{x}$ remains bounded, we have
$$\widetilde{\mathfrak{p}}_{\eps}=\mathcal{O}(1),$$
uniformly with respect to the parameter $0<\eps \leq \eps_0$; while
$$\widetilde{\mathfrak{p}}_{\eps}=\mathcal{O}(|\tilde{x}|^2),$$
when $\abs{\widetilde{x}}\rightarrow +\infty$. It follows from (\ref{eq6.2}) that along the IR-manifold $\Lambda_{\widetilde{\Phi}_{\delta,\eps}}$ the symbol (\ref{eq6.6}) satisfies the following estimate
\begin{equation}\label{eq6.7}\inc
\textrm{Re}\Big(\Big(1-i\tilde{c}\frac{\delta}{|\widetilde{x}^2|}\Big) \widetilde{\mathfrak{p}}_{\eps}\Big(\widetilde{x},\frac{2}{i}\frac{\partial \widetilde{\Phi}_{\delta,\eps}(\widetilde{x})}{\partial \widetilde{x}}\Big)\Big)\geq \frac{\delta}{\tilde{C}}.\num
\end{equation}
in the region where $\abs{\widetilde{x}}\geq 1$; and we recall that associated with the IR-manifold $\Lambda_{\widetilde{\Phi}_{\delta,\eps}}$ is the weighted space
$$H_{\widetilde{\Phi}_{\delta,\eps},\widetilde{h}}(\cc^n),$$
where the index $\tilde{h}$ is here to remind us that the new semiclassical parameter used in the definition of this space is now $\tilde{h}$. Since from (\ref{eq6.4xi}), we have
$$\frac{\widetilde{\Phi}_{\delta,\eps}(\widetilde{x})}{\widetilde{h}} = \frac{\Phi_{\delta,\eps}(x)}{h},$$
we notice that the map
$$u(x)\mapsto \widetilde{u}(\widetilde{x})=\eps^{n/2} u(\sqrt{\eps}\widetilde{x}),$$
then maps unitarily the spaces $H_{\Phi_{\delta,\eps},h}(\cc^n)=H_{\Phi_{\delta,\eps}}(\cc^n)$  to $H_{\widetilde{\Phi}_{\delta,\eps},\widetilde{h}}(\cc^n)$.

The idea is then to make use of the fundamental quantization \textit{vs.} multiplication formula in the space $H_{\widetilde{\Phi}_{\delta,\eps},\widetilde{h}}(\cc^n)$ applied to the operator $\widetilde{P}_{\eps}$. This fundamental formula was established in the analytic case in~\cite{SjDuke} and in the $C^{\infty}$ case in~\cite{HeSjSt} (Section~3.4). Since, as we have already observed, the symbol $\widetilde{\mathfrak{p}}_{\eps}$ may exhibit some quadratic growth, its use requires a bit some care. We verify in the following Proposition that this formula can still be applied in our case.

\bigskip

\begin{prop}\label{prop89}
Let $\psi(\widetilde{x}) \in C^{\infty}_b(\comp^n)$ be such that $\nabla \psi = \mathcal{O}(\abs{\widetilde{x}}^{-1})$ when $|\widetilde{x}| \rightarrow +\infty$. Then the quantization \textit{vs.} multiplication formula holds
\inc\begin{multline*}\label{eq6.8}
(\psi \widetilde{P}_{\eps}\widetilde{u},\widetilde{u})_{\widetilde{\Phi}_{\delta,\eps},\tilde{h}}= \int_{\cc^n} \psi(\widetilde{x}) \widetilde{\mathfrak{p}}_{\eps}\Big(\widetilde{x},\frac{2}{i} \frac{\partial \widetilde{\Phi}_{\delta,\eps}(\widetilde{x})}{\partial \widetilde{x}}\Big) |\widetilde{u}(\widetilde{x})|^2 e^{-\frac{2\widetilde{\Phi}_{\delta,\eps}(\widetilde{x})}{\widetilde{h}}}L(d\widetilde{x}) \\+\mathcal{O}(\widetilde{h})\|\widetilde{u}\|^2_{\widetilde{\Phi}_{\delta,\eps},\tilde{h}}.\num
\end{multline*}
\end{prop}

\bigskip

\begin{proof}
The proof of this proposition follows by inspection of the arguments given in~\cite{HeSjSt} (Section~3.4). Writing a contour integral representation for the quantity $\widetilde{P}_{\eps} \widetilde{u}$ obtained by making a change of scales in the contour integral (\ref{eq3.8}), we shall consider the Taylor expansion of the expression
$$\widetilde{\mathfrak{p}}_{\eps}\Big(\frac{\widetilde{x}+\widetilde{y}}{2},\widetilde{\theta}\Big),$$
along the corresponding contour to get that
\inc\begin{multline*}\label{eq6.9}
\widetilde{\mathfrak{p}}_{\eps}\Big(\frac{\widetilde{x}+\widetilde{y}}{2},\widetilde{\theta}\Big) = \widetilde{\mathfrak{p}}_{\eps}\big(\widetilde{x},\xi(\widetilde{x})\big) + \sum_{j=1}^n f_{j,\eps}(\widetilde{x}) \big(\widetilde{\theta}_j - \xi_j(\widetilde{x})\big)+ \sum_{j=1}^n g_{j,\eps}(\widetilde{x})(\widetilde{y}_j-\widetilde{x}_j) \\ + \widetilde{r}_{\eps}(\widetilde{x},\widetilde{y};\tilde{\theta}), \num
\end{multline*}
where
$$f_{j,\eps}(\widetilde{x}) = \partial_{\theta_j} \widetilde{\mathfrak{p}}_{\eps}\big(\widetilde{x},\xi(\widetilde{x})\big),\ g_{j,\eps}(\widetilde{x}) = \frac{1}{2} \partial_{x_j} \widetilde{\mathfrak{p}}_{\eps}\big(\widetilde{x},\xi(\widetilde{x})\big) \textrm{ and } \xi(\widetilde{x})=\frac{2}{i} \frac{\partial \widetilde{\Phi}_{\delta,\eps}(\widetilde{x})}{\partial \widetilde{x}}.$$
We notice in particular that
\begin{equation}\label{eq6.9.1}\inc
f_{j,\eps}(\widetilde{x})=\mathcal{O}(|\widetilde{x}|), \  \nabla f_{j,\eps}(\widetilde{x})=\mathcal{O}(1),\num
\end{equation}
uniformly with respect to $0<\eps \leq \eps_0$ and $\tilde{x} \in \cc^n$; and that the remainder $\widetilde{r}_{\eps}$ appearing in right hand side of (\ref{eq6.9}) satisfies
$$\widetilde{r}_{\eps}(\widetilde{x},\widetilde{y};\tilde{\theta}) = \mathcal{O}(|\widetilde{x}-\widetilde{y}|^2+h^{\infty}),$$
uniformly with respect to the parameter $0<\eps \leq \eps_0$; as we know that
$$\nabla^2 \widetilde{\mathfrak{p}}_{\eps} \in L^{\infty} \textrm{ and } \widetilde{\theta}-\xi(\widetilde{x})=\mathcal{O}(|\widetilde{x}-\widetilde{y}|),$$
uniformly with respect to the parameter $0<\eps \leq 1$. Arguing as in~\cite{HeSjSt} (Section~3.4), we see that the contribution to the scalar product $(\psi \widetilde{P}_{\eps}\widetilde{u},\widetilde{u})_{\widetilde{\Phi}_{\delta,\eps},\tilde{h}}$ coming from this remainder term $\widetilde{r}_{\eps}$ can be estimated from above by
$$\mathcal{O}(\widetilde{h}) \|\widetilde{u}\|^2_{\widetilde{\Phi}_{\delta,\eps},\tilde{h}}.$$
The contribution from the third term in the right hand side of (\ref{eq6.9}) vanishes as noticed in~\cite{HeSjSt} (Section~3.4). It therefore remains to study the contribution coming from the second term in (\ref{eq6.9}). Here, arguing as in~\cite{HeSjSt} and~\cite{SjDuke}, we see that this term contributes the following sum of integrals
$$\sum_{j=1}^n \int_{\cc^n} \psi(\widetilde{x}) f_{j,\eps}(\widetilde{x}) \big(\widetilde{h}D_{\widetilde{x}_j}-\xi_j(\widetilde{x})\big)\widetilde{u}(\widetilde{x}) \overline{\widetilde{u}(\widetilde{x})} e^{-\frac{2\widetilde{\Phi}_{\delta,\eps}(\widetilde{x})}{\widetilde{h}}}L(d\widetilde{x}).$$
Integrating by parts, as in~\cite{SjDuke} (p.9), and using that $\widetilde{u}$ is a holomorphic function, we find that each of the terms in this sum is equal to
$$-\int_{\cc^n} \widetilde{h} D_{\widetilde{x}_j} \big( \psi(\widetilde{x}) f_{j,\eps} (\widetilde{x})\big)|\widetilde{u}(\widetilde{x})|^2 e^{-\frac{2\widetilde{\Phi}_{\delta,\eps}(\widetilde{x})}{\widetilde{h}}}L(d\widetilde{x}).$$
Considering
$$D_{\widetilde{x}_j}\big(\psi(\widetilde{x}) f_{j,\eps}(\widetilde{x})\big) = D_{\widetilde{x}_j}\psi(\widetilde{x})f_{j,\eps}(\widetilde{x}) + \psi(\widetilde{x})\, D_{\widetilde{x}_j} f_{j,\eps}(\widetilde{x}),$$
we notice that the second term appearing in the right hand side of the last equality is uniformly bounded since from (\ref{eq6.9.1}), we know that $\nabla f_{j,\eps} =\mathcal{O}(1)$. As for the first one, we see using (\ref{eq6.9.1}) and the bound $\nabla \psi = \mathcal{O}(|\widetilde{x}|^{-1})$, that it is bounded as well. This completes the proof of Proposition~\ref{prop89}.
\end{proof}

\bigskip

Considering $\chi$ a $C^{\infty}_b(\comp^n;[0,1])$ function such that $\chi=1$ when $|\widetilde{x}|$ is large enough with
$$\textrm{supp } \chi \subset \{\widetilde{x} \in \cc^n :  |\widetilde{x}|\geq 1\};$$
and
$$m(\widetilde{x})=1-i \tilde{c}\frac{\delta}{|\widetilde{x}|^2},$$
we shall apply Proposition~\ref{prop89} with the function $\psi(\widetilde{x})=\chi(\widetilde{x}) m(\widetilde{x})$ and the operator $\widetilde{P}_{\eps}$ replaced by $\widetilde{P}_{\eps}-\widetilde{h}z$; where the spectral parameter $z\in \comp$ satisfies $\abs{z}\leq C$ for some fixed $C>0$.
From now on, the parameter $\delta>0$ is going to be kept sufficiently small but fixed.
We get from Proposition~\ref{prop89} that the scalar product
\begin{equation}\label{eq6.10}\inc
\big(\chi m (\widetilde{P}_{\eps}-\widetilde{h}z)\widetilde{u},\widetilde{u}\big)_{\widetilde{\Phi}_{\delta,\eps},\tilde{h}}, \num
\end{equation}
is equal to
$$\int_{\cc^n} \chi(\widetilde{x})m(\widetilde{x})\widetilde{\mathfrak{p}}_{\eps}\Big(\widetilde{x},\frac{2}{i} \frac{\partial \widetilde{\Phi}_{\delta,\eps}(\widetilde{x})}{\partial \widetilde{x}}\Big) \abs{\widetilde{u}(\widetilde{x})}^2e^{-\frac{2\widetilde{\Phi}_{\delta,\eps}(\widetilde{x})}{\widetilde{h}}}L(d\widetilde{x}) +\mathcal{O}(\widetilde{h})\|\widetilde{u}\|^2_{\widetilde{\Phi}_{\delta,\eps},\tilde{h}}.$$
Since
\begin{multline*}
\textrm{Re}\big(\chi m (\widetilde{P}_{\eps}-\widetilde{h}z)\widetilde{u},\widetilde{u}\big)_{\widetilde{\Phi}_{\delta,\eps},\tilde{h}} \\
= \int_{\cc^n} \chi(\widetilde{x}) \textrm{Re}\Big[m(\widetilde{x})\widetilde{\mathfrak{p}}_{\eps}\Big(\widetilde{x},\frac{2}{i} \frac{\partial \widetilde{\Phi}_{\delta,\eps}(\widetilde{x})}{\partial \widetilde{x}}\Big)\Big]|\widetilde{u}(\widetilde{x})|^2 e^{-\frac{2\widetilde{\Phi}_{\delta,\eps}(\widetilde{x})}{\widetilde{h}}}L(d\widetilde{x}) +\mathcal{O}(\widetilde{h})\|\widetilde{u}\|^2_{\widetilde{\Phi}_{\delta,\eps},\tilde{h}},
\end{multline*}
we deduce from the Cauchy-Schwarz inequality and the estimate (\ref{eq6.7}) holding on the support of the function $\chi$ that
$$\int_{\cc^n} \chi(\widetilde{x}) |\widetilde{u}(\widetilde{x})|^2 e^{-\frac{2\widetilde{\Phi}_{\delta,\eps}(\widetilde{x})}{\widetilde{h}}}L(d\widetilde{x}) \leq \mathcal{O}(1) \|(\widetilde{P}_{\eps}-\widetilde{h}z)\widetilde{u}\|_{\widetilde{\Phi}_{\delta,\eps},\tilde{h}} \|\widetilde{u}\|_{\widetilde{\Phi}_{\delta,\eps},\tilde{h}} + \mathcal{O}(\widetilde{h})\|\widetilde{u}\|^2_{\widetilde{\Phi}_{\delta,\eps},\tilde{h}}.$$
We can now come back to the original variables $x=\sqrt{\eps}\widetilde{x}$ and obtain that
\begin{multline*}
\eps \int_{\cc^n} \chi\Big(\frac{x}{\sqrt{\eps}}\Big) |u(x)|^2 e^{-\frac{2\Phi_{\delta,\eps}(x)}{h}}L(dx) \leq \mathcal{O}(1) \|(\mathfrak{p}_0^w(x,hD_x)-hz)u\|_{\Phi_{\delta,\eps}}\|u\|_{\Phi_{\delta,\eps}}\\ +\mathcal{O}(h)\|u\|^2_{\Phi_{\delta,\eps}}.
\end{multline*}
It follows from (\ref{xi13}) that the operator $P_0=P_0^w(x,hD_x;h)$ also fulfills
$$\eps \int_{\cc^n} \chi\Big(\frac{x}{\sqrt{\eps}}\Big) |u(x)|^2 e^{-\frac{2\Phi_{\delta,\eps}(x)}{h}}L(dx) \leq \mathcal{O}(1) \|(P_0-hz)u\|_{\Phi_{\delta,\eps}}\|u\|_{\Phi_{\delta,\eps}} +\mathcal{O}(h)\|u\|^2_{\Phi_{\delta,\eps}}.$$
Finally, by recalling that $\eps=Ah$, we can rewrite this estimate as
\inc\begin{multline}\label{eq6.11}
h \int_{\cc^n} \chi\Big(\frac{x}{\sqrt{A h}}\Big)|u(x)|^2 e^{-\frac{2\Phi_{\delta,\eps}(x)}{h}}L(dx)  \leq \mathcal{O}(1) \|(P_0-hz)u\|_{\Phi_{\delta,\eps}}\|u\|_{\Phi_{\delta,\eps}}
 \\ + \mathcal{O}\Big(\frac{h}{A}\Big)\|u\|_{\Phi_{\delta,\eps}}^2, \num
\end{multline}
where $A \gg 1$ is a large parameter remaining to be chosen.

\section{Proof of Theorem~\ref{theo}}\label{proof}
\init

In this section, we shall glue together the local resolvent estimates that we established first in a tiny neighborhood of the doubly characteristic set (Section~\ref{tiny}) and then in the exterior region considered in the previous section in order to complete the proof of Theorem~\ref{theo}.

By considering the same cutoff function $\chi_0 \in C_0^{\infty}(\comp^n,[0,1])$ as in Section~\ref{tiny} and going back to (\ref{k41}) with $k=1/2$, we deduce from this estimate and the triangle inequality that there exists $C>0$ such that
\inc\begin{align*}\label{chel1}
& \ \Big\| (h+d^2)^{\frac{1}{2}}\chi_0\Big(\frac{x}{\sqrt{Ah}}\Big) u\Big\|_{\Phi_{\delta,\eps}} \num\\
\leq & \ C \Big\| (h+d^2)^{-\frac{1}{2}}\chi_0\Big(\frac{x}{\sqrt{Ah}}\Big)  \big(\mathfrak{p}_0^w(x,hD_x)+h\mathfrak{p}_1^w(x,hD_x)-hz\big)u\Big\|_{\Phi_{\delta,\eps}}\\
+ & \ C \Big\| (h+d^2)^{-\frac{1}{2}}\chi_0\Big(\frac{x}{\sqrt{Ah}}\Big)  \big(\mathfrak{p}_0^w(x,hD_x)+h\mathfrak{p}_1^w(x,hD_x)-Q_0-hp_1(0,0)\big)u\Big\|_{\Phi_{\delta,\eps}}\\
+ & \ \frac{\tilde{C}}{\sqrt{A}} \Big\|(h+d^2)^{\frac{1}{2}}\un_{K}\Big(\frac{x}{\sqrt{Ah}}\Big)u\Big\|_{\Phi_{\delta,\eps}},
\end{align*}
for all $z$ varying in a compact set such that the set $z-p_1(0,0)$ does not contain any eigenvalues of the operator $Q_0|_{h=1}$ described in (\ref{k27}). We recall that $Q_0$ stands here for the operator appearing in (\ref{k36}). Since on the supports of the functions $\chi_0(x/\sqrt{Ah})$ and $\un_{K}(x/\sqrt{Ah})$, one can estimate the quantity $h+d^2$ by
$$h \leq h+d^2 \leq c_1 Ah,$$
with $c_1>0$ a positive constant independent of the parameters $A$ and $h$; we deduce from (\ref{k36}) and (\ref{chel1}) that there exist some positive constants $c_{2}$ and $c_{3,A}$, where the constant $c_2$  is independent of the parameters $A$ and $h$, whereas $c_{3,A}$ may actually depend on $A$ but not on $h$; such that
\inc\begin{multline*}\label{k51xi}
h^{\frac{1}{2}}\Big\| \chi_0\Big(\frac{x}{\sqrt{Ah}}\Big) u\Big\|_{\Phi_{\delta,\eps}} \leq c_2 h^{-\frac{1}{2}} \big\|\big(\mathfrak{p}_0^w(x,hD_x)+h\mathfrak{p}_1^w(x,hD_x)-hz\big)u\big\|_{\Phi_{\delta,\eps}} \\ +c_{3,A} h \|u\|_{\Phi_{\delta,\eps}}  +c_2 h^{\frac{1}{2}}  \Big\|\un_{K}\Big(\frac{x}{\sqrt{Ah}}\Big)u\Big\|_{\Phi_{\delta,\eps}}. \num
\end{multline*}
It follows from (\ref{xi13}) that the operator $P_0=P_0^w(x,hD_x;h)$ also fulfills
\inc\begin{multline*}\label{k51}
h^{\frac{1}{2}}\Big\| \chi_0\Big(\frac{x}{\sqrt{Ah}}\Big) u\Big\|_{\Phi_{\delta,\eps}} \leq c_2 h^{-\frac{1}{2}} \big\|\big(P_0-hz\big)u\big\|_{\Phi_{\delta,\eps}} +c_{3,A} h \|u\|_{\Phi_{\delta,\eps}}  \\ +c_2 h^{\frac{1}{2}}  \Big\|\un_{K}\Big(\frac{x}{\sqrt{Ah}}\Big)u\Big\|_{\Phi_{\delta,\eps}}, \num
\end{multline*}
where $c_{3,A}$ stands for a new constant depending on the parameter $A$ but not on $h$.
Recalling that here $K$ stands for a fixed neighborhood of the support of the function $\nabla \chi_0$, we get from (\ref{k51}) upon squaring that
\inc\begin{multline}\label{eq51.1}
h \Big\| \chi_0\Big(\frac{x}{\sqrt{Ah}}\Big) u\Big\|_{\Phi_{\delta,\eps}}^2
\leq \frac{\mathcal{O}(1)}{h} \|(P_0-hz)u\|_{\Phi_{\delta,\eps}}^2 + \mathcal{O}_A(h^2)\|u\|_{\Phi_{\delta,\eps}}^2 \\ + \mathcal{O}(h) \Big\|\un_{K}\Big(\frac{x}{\sqrt{Ah}}\Big)u\Big\|_{\Phi_{\delta,\eps}}^2. \num
\end{multline}
On the other hand, we get from the estimate (\ref{eq6.11}) that
\inc\begin{multline}\label{eq51.2}
h \int_{\cc^n} \chi\Big(\frac{x}{\sqrt{Ah}}\Big) |u(x)|^2 e^{-\frac{2\Phi_{\delta,\eps}(x)}{h}}L(dx)+ h \Big\|\un_{K}\Big(\frac{x}{\sqrt{Ah}}\Big)u\Big\|_{\Phi_{\delta,\eps}}^2 \\
\leq \mathcal{O}(1) \|(P_0-hz)u\|_{\Phi_{\delta,\eps}} \|u\|_{\Phi_{\delta,\eps}}+\mathcal{O}\Big(\frac{h}{A}\Big)\|u\|_{\Phi_{\delta,\eps}}^2, \num
\end{multline}
and we then obtain by collecting (\ref{eq51.1}) and (\ref{eq51.2}) that
\begin{multline*}
h \int_{\cc^n}|u(x)|^2 e^{-\frac{2\Phi_{\delta,\eps}(x)}{h}}L(dx) \leq \mathcal{O}(1) \|(P_0-hz)u\|_{\Phi_{\delta,\eps}}\|u\|_{\Phi_{\delta,\eps}} \\
 +\frac{\mathcal{O}(1)}{h} \|(P_0-hz)u\|_{\Phi_{\delta,\eps}}^2 +  \Big(\mathcal{O}\Big(\frac{h}{A}\Big)+\mathcal{O}_A(h^2)\Big)\|u\|_{\Phi_{\delta,\eps}}^2.
\end{multline*}
Here we have also used that we arrange, as we may, that $\chi+\chi_0^2\geq 1$ on $\cc^n$.
By multiplying by the parameter $h$ and using that
$$\mathcal{O}(1) h \|(P_0-hz)u\|_{\Phi_{\delta,\eps}} \|u\|_{\Phi_{\delta,\eps}} \leq \mathcal{O}(1) \|(P_0-hz)u\|_{\Phi_{\delta,\eps}}^2 + \frac{h^2}{2}\|u\|_{\Phi_{\delta,\eps}}^2,$$
we get that
\begin{multline*}
h^2 \int_{\cc^n} |u(x)|^2 e^{-\frac{2\Phi_{\delta,\eps}(x)}{h}}L(dx) \leq \mathcal{O}(1) \|(P_0-hz)u\|_{\Phi_{\delta,\eps}}^2 + \frac{h^2}{2}\|u\|_{\Phi_{\delta,\eps}}^2 \\
+ \Big(\mathcal{O}\Big(\frac{h^2}{A}\Big) + \mathcal{O}_A(h^3)\Big)\|u\|_{\Phi_{\delta,\eps}}^2.
\end{multline*}
By now choosing the parameter $A$ sufficiently large, but fixed; and then considering a positive constant $h_0>0$ sufficiently small, $0<h_0 \ll 1$, depending on the choice done for the value of the constant $A$, we obtain that for all $0 < h \leq h_0$,
\inc\begin{equation}\label{eq51.3}
h \|u\|_{\Phi_{\delta,\eps}} \leq \mathcal{O}(1)\|(P_0-hz)u\|_{\Phi_{\delta,\eps}}. \num
\end{equation}
Let us underline that one can then replace in this estimate the weight function $\Phi_{\delta,\eps}$ by the standard quadratic weight $\Phi_0$ defined in (\ref{eq3.2.1b}) at the expense of an $\mathcal{O}(1)$--loss. Indeed, this follows from the fact that according to (\ref{ju6}), the two associated $L^2$-norms are equivalent since we have
$$\textrm{exp}\big(-\mathcal{O}(A)\big) \leq e^{-\frac{2\Phi_{\delta,\eps}}{h}} e^{\frac{2\Phi_0}{h}}\leq \textrm{exp}\big(\mathcal{O}(A)\big).$$
It is now easy to complete the proof of Theorem~\ref{theo} . In doing so, it is sufficient to go back to the $L^2$--side by undoing the FBI transform $T$ in (\ref{eq51.3}) to get that for all $u\in L^2(\rr^n)$ and $0<h \leq h_0$,
\begin{equation}\label{eq51.5xi}\inc
h\|u\|_{L^2} \leq \mathcal{O}(1) \|(P-hz)u\|_{L^2}, \num
\end{equation}
for all $z$ varying in a compact set such that the set $z-p_1(0,0)$ does not contain any eigenvalues of the operator $q_1^w(x,D_x)$ described in (\ref{k27}). This ends the proof of Theorem~\ref{theo} since we recall that we were considering here, for simplicity only, the case when the doubly characteristic set is reduced to an unique point $X_1=(0,0)$.
In the general case when the doubly characteristic set is composed of a finite number of points
$$p_0^{-1}(0)=\{X_1,...,X_N\},$$
a simple adaptation of the previous arguments allows to establish similar estimates as (\ref{eq51.1}) near each doubly characteristic point $X_j$ when the spectral parameter $z$ stays in a compact set as described in the statement of Theorem~\ref{theo}. We then conclude the proof of Theorem~\ref{theo} by using a similar a priori estimate as (\ref{eq51.2}) in the exterior region.~$\Box$

\end{document}